\documentclass[11pt]{article}

\usepackage{amsfonts} 
\usepackage{amsmath,hhline,latexsym,mathrsfs}
\usepackage{amssymb}
\usepackage{relsize}
\DeclareMathAlphabet\mathbfcal{OMS}{cmsy}{b}{n}
\usepackage{mathtools}

\usepackage[latin1]{inputenc}

\usepackage{hyperref}

\usepackage{color}
\usepackage{graphicx}
\usepackage{comment}

\newtheorem{theorem}{Theorem}
\newtheorem{prop}{Proposition}
\newtheorem{lemma}{Lemma}

\newtheorem{remark}{Remark}

\newcommand{\dis}{\displaystyle}
\newcommand{\R}{\mathbb{R}}
\newcommand{\N}{\mathbb{N}}
\newcommand{\A}{\mathcal{A}}

\begin{document}

\title{Approximation of null controls for semilinear heat equations using a least-squares approach}

\author{J\'er\^ome Lemoine \thanks{Laboratoire de math\'ematiques Blaise Pascal, Universit\'e Clermont Auvergne, UMR CNRS 6620, Campus des C\'ezeaux, 3, place Vasarely, 63178 Aubi\`ere, France. e-mail: jerome.lemoine@uca.fr.}, 
\ \ Irene Mar\'{\i}n-Gayte\thanks{Departamento EDAN, Universidad de Sevilla, Campus Reina Mercedes, 41012, Sevilla, Spain. e-mail: imgayte@us.es.},
\ \ Arnaud M\"unch\thanks{Laboratoire de math\'ematiques Blaise Pascal, Universit\'e Clermont Auvergne, UMR CNRS 6620, Campus des C\'ezeaux, 3, place Vasarely, 63178 Aubi\`ere, France. e-mail: arnaud.munch@uca.fr (Corresponding author).}
 }

\maketitle


\vspace*{0.5cm}

\begin{abstract}
The null distributed controllability of the semilinear heat equation $y_t-\Delta y + g(y)=f \,1_{\omega}$, assuming that $g$ satisfies the growth condition $g(s)/(\vert s\vert \log^{3/2}(1+\vert s\vert))\rightarrow 0$ as $\vert s\vert \rightarrow \infty$ and that $g^\prime\in L^\infty_{loc}(\mathbb{R})$ has been obtained by Fern\'andez-Cara and Zuazua in 2000. The proof based on a fixed point argument makes use  of precise estimates of the observability constant for a linearized heat equation. It does not provide however an explicit construction of a null control.   
Assuming that $g^\prime\in W^{s,\infty}(\mathbb{R})$ for one $s\in (0,1]$, we construct an explicit sequence converging strongly to a null control for the solution of the semilinear equation. The method, based on a least-squares approach, generalizes Newton type methods and guarantees the convergence whatever be the initial element of the sequence. In particular, after a finite number of iterations, the convergence is super linear with a rate equal to $1+s$. Numerical experiments in the one dimensional setting support our analysis. 
 \end{abstract}
\vspace*{0,5in}
 


\section{Introduction}

Let $\Omega\subset \mathbb{R}^d$, $1\leq d\leq 3$, be a bounded connected open set whose boundary $\partial\Omega$ is Lipschitz. Let $\omega$ be any non-empty open set of $\Omega$ and let $T>0$. We note $Q_T=\Omega\times (0,T)$, $q_T=\omega\times (0,T)$ and $\Sigma_T=\partial\Omega\times (0,T)$.  We are concerned with the null controllability problem for the following semilinear heat equation
\begin{equation}
\label{heat-NL}
\left\{
\begin{aligned}
& y_t - \Delta y +  g(y)= f 1_{\omega} \quad  \textrm{in}\quad Q_T,\\
& y=0 \,\, \textrm{on}\,\, \Sigma_T, \quad y(\cdot,0)=u_0 \,\, \textrm{in}\,\, \Omega,
\end{aligned}
\right.
\end{equation}
where $u_0\in L^2(\Omega)$ is the initial state of $y$ and $f\in L^2(q_T)$ is a {\it control} function. We assume moreover that the nonlinear function $g:\mathbb{R} \mapsto \mathbb{R}$ is, at~least, locally Lipschitz-continuous. Following \cite{EFC-EZ}, we will also assume for simplicity that $g$ satisfies
   \begin{equation}
   \vert  g^{\prime}(s) \vert \leq C (1+\vert s\vert^m) \quad \textrm{a.e., with}\,\,\, 1 \leq m \leq 1+4/d. \label{cond_f_5}
   \end{equation}
   Under this condition, (\ref{heat-NL}) possesses exactly one local in time solution.
   Moreover, under the growth condition
   \begin{equation}
\vert g(s)\vert \leq C(1+\vert s\vert \log(1+\vert s\vert)) \quad \forall s\in \mathbb{R}, \label{growth-f}
   \end{equation}
   the solutions to (\ref{heat-NL}) are globally defined in $[0,T]$ and one has
   \begin{equation}\label{globalT}
   y \in C^0([0,T]; L^2(\Omega))\cap L^2(0,T; H_0^1(\Omega)),
   \end{equation}
 see~\cite{CazenaveHaraux}.
   Recall that, without a growth condition of the kind (\ref{growth-f}), the solutions to~(\ref{heat-NL}) can blow up before $t=T$;
   in general, the blow-up time depends on~$g$ and the size of $\Vert u_0\Vert_{L^2(\Omega)}$.

  The system (\ref{heat-NL}) is said to be {\it controllable} at time $T$ if, for any $u_0\in L^2(\Omega)$ and any globally defined bounded trajectory $y^{\star} \in C^0([0,T]; L^2(\Omega))$
   (corresponding to the data $u_0^\star \in L^2(\Omega)$ and $f^{\star}\in L^2(q_T)$), there exist controls $f \in L^2(q_T)$ and associated states~$y$ that are again globally defined in $[0,T]$ and satisfy~\eqref{globalT} and
   \begin{equation}\label{stateT}
   y(x,T) = y^{\star}(x,T), \quad x \in \Omega.
   \end{equation}
   We refer to \cite{coron-book} for an overview of control problems in nonlinear situations. The uniform controllability strongly depends on the nonlinearity $g$. Fern\'andez-Cara and Zuazua proved in \cite{EFC-EZ} that if $g$ is too ``super-linear" at infinity, then, for some initial data, the control cannot compensate the blow-up phenomenon occurring in $\Omega \backslash \overline{\omega}$:

  
\begin{theorem}[\cite{EFC-EZ}]\label{nullcontrollheatmoins}
   There exist locally Lipschitz-continuous functions $g$ with $g(0)=0$ and
   \begin{equation}
\vert g(s)\vert \sim \vert s\vert \log^p (1+\vert s \vert) \quad \textrm{as}\quad \vert s\vert \rightarrow \infty,
\quad p > 2, \nonumber
   \end{equation}
such that $(\ref{heat-NL})$ fails to be controllable for all $T>0$.
\end{theorem}
On the other hand, Fern\'andez-Cara and Zuazua also proved that if $p$ is small enough, then the controllability holds true uniformly.   
\begin{theorem}[\cite{EFC-EZ}]\label{nullcontrollheatplus}
   Let $T>0$ be given.
   Assume that $(\ref{heat-NL})$ admits at least one solution $y^{\star}$, globally defined in $[0,T]$ and bounded in $Q_T$.
   Assume that $g:\mathbb{R} \mapsto \mathbb{R}$ is locally Lipschitz-continuous and satisfies $(\ref{cond_f_5})$ and
   \begin{equation} \label{asymptotic_g}
\frac{g(s)}{\vert s\vert \log^{3/2}(1+\vert s\vert)}\rightarrow 0 \quad \textrm{as}\quad \vert s\vert \rightarrow \infty.
   \end{equation}
   Then $(\ref{heat-NL})$ is controllable at time $T$.
\end{theorem}
Therefore, if $|g(s)|$ does not grow at infinity faster than $\vert s \vert \log^{p}(1+\vert s\vert)$ for any $p<3/2$, then (\ref{heat-NL}) is controllable. This result extends \cite{EFC97} obtaining the uniform controllability for any  $p<1$. We also mention \cite{VB1} which gives the same result assuming additional sign condition on $g$, namely $g(s)s\geq -C (1+s^2)$ for all $s\in \mathbb{R}$ and some $C>0$. The problem remains open when $g$ behaves at infinity like $\vert s\vert \log^p(1+\vert s\vert)$ with $3/2\leq p\leq 2$. We mention however the recent work of LeBalc'h \cite{KLB} where uniform controllability results are obtained for $p\leq 2$ assuming additional sign conditions on $g$, notably that $g(s)>0$ for $s>0$ and $g(s)<0$ for $s<0$. This condition is not satisfied for $g(s)=-s\, \log^p(1+\vert s\vert)$. 
 Let us also mention \cite{coron-trelat} in the context of Theorem \ref{nullcontrollheatmoins} where a positive boundary controllability result is proved for a specific class of initial and final data and $T$ large enough.

In the sequel, for simplicity, we shall assume that $g(0)=0$ and that $f^\star\equiv 0, u_0^\star\equiv 0$ so that $y^\star$ is the null trajectory. The proof given in \cite{EFC-EZ} is based on a fixed point method. Precisely, it is shown that the operator $\Lambda:L^\infty(Q_T)\to L^\infty(Q_T)$, where $y_z:=\Lambda z$ is a null controlled solution of the linear boundary value problem 
\begin{equation}
\label{NL_z}
\left\{
\begin{aligned}
& y_{z,t} - \Delta y_z +  y_z \,\tilde{g}(z)= f_z 1_{\omega} \quad  \textrm{in}\quad Q_T\\
& y_z=0 \,\, \textrm{on}\,\, \Sigma_T, \quad y_z(\cdot,0)=u_0 \quad \textrm{in}\quad \Omega
\end{aligned}
\right., 
\qquad
\tilde{g}(s):=
\left\{ 
\begin{aligned}
& g(s)/s & s\neq 0,\\
& g^{\prime}(0) & s=0,
\end{aligned}
\right.
\end{equation}
maps a closed ball $B(0,M) \subset L^\infty(Q_T)$ into itself, for some $M>0$. The Kakutani's theorem then provides the existence of at least one fixed point for the operator $\Lambda$, which is also a controlled solution for (\ref{heat-NL}).

The main goal of this work is to determine an approximation of the controllability problem associated to~(\ref{heat-NL}), that is to construct an explicit sequence $(f_k)_{k\in \mathbb{N}}$ converging strongly toward a null control for~(\ref{heat-NL}). A natural strategy is to take advantage of the method used in \cite{KLB, EFC-EZ} and consider the Picard iterates associated with the operator $\Lambda$: $y_{k+1}=\Lambda(y_k)$, $k\geq 0$ initialized with any element $y_0\in B(0,M)$. The sequence of controls is then $(f_k)_{k\in \mathbb{N}}$ so that $f_k\in L^2(q_T)$ is a null control for $y_k$ solution of 
\begin{equation}
\label{NL_z_k}
\left\{
\begin{aligned}
& y_{k,t} - \Delta y_{k} +  y_{k} \,\tilde{g}(y_{k-1})= f_{k} 1_{\omega} \quad  \textrm{in}\quad Q_T,\\
& y_{k}=0 \,\, \textrm{on}\,\, \Sigma_T, \quad y_{k}(\cdot,0)=u_0 \,\, \textrm{in}\,\, \Omega.
\end{aligned}
\right. 
\end{equation}
Numerical experiments for $d=1$ reported in \cite{EFC-AM} exhibit the non convergence of the sequences $(y_k)_{k\in\mathbb{N}}$ and  $(f_k)_{k\in\mathbb{N}}$ for some initial conditions large enough. This phenomenon is related to the fact that the operator $\Lambda$ is \textit{a priori} not contractant. We also refer to \cite{BoyerCanum2012} where this strategy is implemented. 
Still in the one dimensional case, a least-squares type approach, based on the minimization over $L^2(Q_T)$ of the functional $R:L^2(Q_T)\to \mathbb{R}^+$ defined by $R(z):=\Vert z-\Lambda(z)\Vert_{L^2(Q_T)}$
is introduced and analyzed in \cite{EFC-AM}. Assuming that $\tilde{g}\in C^1(\mathbb{R})$ and $g^{\prime}\in L^\infty(\mathbb{R})$, it is proved first that $R\in C^1(L^2(Q_T);\mathbb{R}^+)$ and secondly that, if $\Vert u_0\Vert_{L^\infty(\Omega)}$ is small enough, then any critical point for $R$ is a fixed point for $\Lambda$. Under this smallness  assumption on the data, numerical experiments reported in \cite{EFC-AM} display the convergence of minimizing sequences for $R$ (based on a gradient method) and a better behavior than the Picard iterates. The analysis of convergence is however not performed.  As is usual for nonlinear problems and considered in \cite{EFC-AM}, we may also employ a Newton type method to find a zero of the mapping $\widetilde{F}: Y \mapsto W$ defined by
   \begin{equation}\label{def-F}
\widetilde{F}(y,f) = (y_t - \Delta y  + g(y) - f 1_{\omega}, y(\cdot\,,0) - u_0, y(\cdot,T)) \quad \forall (y,f) \in Y
   \end{equation}
 for some appropriate Hilbert spaces $Y$ and $W$ (see below).
%
It is shown for $d=1$ in \cite{EFC-AM} that, if $g\in C^1(\mathbb{R})$ and $g^{\prime}\in L^\infty(\mathbb{R})$, then $\widetilde{F}\in C^1(Y;W)$ allowing to derive the Newton iterative sequence: given $(y_0,f_0)$ in $Y$, define the sequence $(y_k,f_k)_{k\in \mathbb{N}}$ iteratively as follows $(y_{k+1},f_{k+1})=(y_k,f_k)-(Y_k,F_k)$ where $F_k$ is a control for $Y_k$ solution of  
\begin{equation}
\label{Newton-nn}
\left\{
\begin{aligned}
&Y_{k,t}-\Delta Y_{k} +  g^{\prime}(y_k)\,Y_{k} = F_{k}\, 1_{\omega} + y_{k,t} - \Delta y_k  + g(y_k) - f_k 1_{\omega},
                                                                                    & \quad  \textrm{in}\quad Q_T,\\
& Y_k=0,                                                              & \quad  \textrm{on}\quad \Sigma_T,\\ 
& Y_k(\cdot,0)=u_0-y_k(\cdot,0),                                     & \quad  \textrm{in}\quad \Omega
\end{aligned}
\right.
   \end{equation}
   such that $Y_k(\cdot,T)=-y_k(\cdot,T)$ in $\Omega$. Once again, numerical experiments for $d=1$ in \cite{EFC-AM} exhibits the lack of convergence of the Newton method for large enough initial condition, for which the solution $y$ is not close enough to the zero trajectory. As far as we know, the construction of a convergent approximation $(f_{k})_{k\in \mathbb{N}}$ in the general case where the initial data to be controlled is arbitrary in $L^2(\Omega)$ remains an open issue. 
  Still assuming that $g^{\prime}\in L^\infty(\mathbb{R})$ and in addition that there exists one $s$ in $(0,1]$ such that $\sup_{a,b\in \mathbb{R}, a\neq b} \frac{\vert g^\prime(a)-g^\prime(b)\vert}{\vert a-b\vert^s}< \infty$, we construct, for any initial data $u_0\in L^2(\Omega)$, a strongly convergent sequence $(f_k)_{k\in \mathbb{N}}$ toward a control for (\ref{heat-NL}). Moreover, after a finite number of iterates related to the norm $\Vert g^{\prime}\Vert_{L^\infty(\mathbb{R})}$, the convergence is super linear with a rate equal to $1+s$. This is done (following and improving \cite{AM-PP-2014} devoted to a linear case) by introducing a quadratic functional which measures how a pair $(y,f) \in Y$ is close to a controlled solution for (\ref{heat-NL}) and then by determining a particular minimizing sequence enjoying the announced property. A natural example of so-called error (or least-squares) functional is given by $\widetilde{E}(y,f):=\frac{1}{2}\Vert \widetilde{F}(y,f)\Vert^2_W$ to be minimized over $Y$. In view of controllability results for (\ref{heat-NL}), the non-negative functional $\widetilde{E}$ achieves its global minimum equal to zero for any control pair $(y,f)\in Y$ of (\ref{heat-NL}).

  The paper is organized as follows. In Section \ref{sec_control_linearized}, we first derive a controllability result for a linearized wave equation with potential in $L^\infty(Q_T)$ and source term in $L^2(0,T;H^{-1}(\Omega))$. Then, in Section \ref{sec_LS}, we define the least-squares functional $E$ and the corresponding optimization problem \eqref{extremal_problem} over the Hilbert space $\mathcal{A}$. We show that $E$ is Gateaux-differentiable over $\mathcal{A}$ and that any critical point $(y,f)$ for $E$ for which $g^\prime(y)$ belongs to $L^\infty(Q_T)$ is also a zero of $E$ (see Proposition \ref{proposition3}). This is done by introducing a descent direction $(Y^1,F^1)$ for $E(y,f)$ for which $E^\prime(y,f)\cdot (Y^1,F^1)$ is proportional to $E(y,f)$. Then, assuming that the nonlinear function $g$ is such that $g^\prime$ belongs to $W^{s,\infty}(\mathbb{R})$ for one $s$ in $(0,1]$, we determine a minimizing sequence based on $(Y^1,F^1)$ which converges strongly to a controlled pair for the semilinear heat equation (\ref{heat-NL}). Moreover, we prove that after a finite number of iterates, the convergence enjoys a rate equal to $1+s$ (see Theorem \ref{ths1} for $s=1$ and Theorem \ref{ths} for $s\in (0,1)$). We also emphasize that this least-squares approach coincides with the damped Newton method one may use to find a zero of a mapping similar to $\widetilde{F}$ mentioned above; we refer to Remark \ref{link_DampedNewtonMethod}. This explains the convergence of our approach with a super-linear rate. Section \ref{sec_numeric} gives some numerical illustrations of our result in the one dimensional case and a nonlinear function $g$ for which $g^\prime\in W^{1,\infty}(\mathbb{R})$.  We conclude in Section \ref{conclusion} with some perspectives. 
  As far as we know, the analysis of convergence presented in this work, though some restrictive hypotheses on the nonlinear function $g$, is the first one in the context of controllability for partial differential equations.

Along the text,  we shall denote by $\Vert \cdot \Vert_{\infty} $ the usual norm in $L^\infty(\mathbb{R})$, $(\cdot,\cdot)_X$ the scalar product of $X$ (if $X$ is a Hilbert space) and by $\langle \cdot, \cdot \rangle_{X,Y}$ the duality product between the spaces $X$ and $Y$.

\section{A controllability result for a linearized heat equation with $L^2(H^{-1})$ right hand side}\label{sec_control_linearized}

We give in this section a controllability result for a linear heat equation with potential in $L^\infty(Q_T)$ and right hand side in $L^2(0,T;H^{-1}(\Omega))$.
As this work concerns the null controllability of parabolic equation, we shall make use  of Carleman type weights introduced in this context notably in \cite{fursikov-imanuvilov} (we also refer to \cite{EFC-SG} for a review). Here, we assume that such weights $\rho$, $\rho_0$, $\rho_1$ and $\rho_2$ blow up as $t \to T^-$ and satisfy:
   \begin{equation}
\label{Hyp-rho-1}
\left\{
\begin{array}{l}
\dis \hbox{$\rho = \rho(x,t)$, $\rho_{0} = \rho_{0}(x,t)$, $\rho_1 = \rho_1(x,t)$ and $\rho_2 = \rho_2(x,t)$   are continuous and $\geq \rho_{*} > 0$~in $Q_{T}$ } \\
\rho, \rho_{0},\rho_1,\rho_2 \in L^\infty(Q_{T-\delta}) \quad \forall \delta > 0.
\end{array}
\right.
   \end{equation}
   
   Precisely, we will take $\rho_0=(T-t)^{3/2} \rho$, $\rho_1=(T-t)\rho$   and $\rho_2=(T-t)^{1/2}\rho$ where $\rho$ is defined as follow
  \begin{equation}\label{def_weight}
  \rho(x,t)=\exp\Big(\frac{s\beta(x)}{\ell(t)}\Big),\quad s\ge C(\Omega,\omega,T,\|g^\prime\|_\infty)
  \end{equation}
   with $ 
   \ell(t)=\begin{cases}
   t(T-t)& \mbox{ si }t\ge T/4\\ 
   3T^2/16& \mbox{ si } 0\le t<T/4\end{cases}
   $.
Here $\beta(x)=\exp(2\lambda m\|\eta^0\|_\infty)-\exp(\lambda(m\|\eta^0\|_\infty+\eta^0(x)))$, $m>1$, $\eta^0\in \mathcal{C}(\overline{\Omega})$ satisfies $\eta^0>0$ in $\Omega$, $\eta^0=0$ on $\partial\Omega$ and $|\nabla \eta^0|>0$ in $\overline{\Omega\backslash\omega}$ (see \cite{EFC-SG}, Lemma 1.2, p.1401).

In the next section, we shall make use the following controllability result.

\begin{prop}\label{controllability_result}
Assume $A\in L^{\infty}(Q_T)$,  $\rho_2B\in L^2(0,T;H^{-1}(\Omega))$ and $z_0\in L^2(\Omega)$. Then there exists a control $v\in L^2(\rho_0,q_T)$ such that the weak solution $z$ of 
\begin{equation}
\label{heat_z}
\left\{
\begin{aligned}
& z_t - \Delta z +  A z  = v 1_{\omega} + B \quad \textrm{in}\quad Q_T,\\
& z=0 \,\, \textrm{on}\,\, \Sigma_T, \quad z(\cdot,0)=z_0 \,\, \textrm{in}\,\, \Omega
\end{aligned}
\right.
\end{equation}
satisfies
\begin{equation}
\label{heat_z1}
z(\cdot,T)=0  \hbox{ in }\Omega.
\end{equation}
  Moreover, the unique control $u$ which minimizes together with the corresponding solution $z$ the functional 
$J:L^2(\rho,Q_T)\times L^2(\rho_0,q_T)\to \mathbb{R}^+$ defined by $J(z,v):=\frac12\Vert \rho\,  z\Vert^2_{L^2(Q_T)} + \frac12\Vert \rho_0\,  v\Vert^2_{L^2(q_T)}$
satisfies the following  estimate
\begin{equation}\label{estimation1}
\Vert \rho\, z\Vert_{L^2(Q_T)}+ \Vert \rho_0\,v\Vert_{L^2(q_T)} \leq C \biggl(\Vert \rho_2 B\Vert_{L^2(0,T;H^{-1}(\Omega))} + \Vert z_0\Vert_{L^2(\Omega)}\biggr) 
\end{equation}
for some  constant $C=C(\Omega,\omega,T, \Vert A\Vert_{\infty})$.

The controlled solution also satisfies, for some  constant $C=C(\Omega,\omega,T, \Vert A\Vert_{\infty})$, the estimate 
\begin{equation}\label{estimation2}
\Vert \rho_1 z\Vert_{L^\infty(0,T;L^2(\Omega))}+ \Vert \rho_1\nabla z\Vert_{L^2(Q_T)^d}\leq C \biggl(\Vert \rho_2\,B\Vert_{L^2(0,T;H^{-1}(\Omega))} + \Vert z_0\Vert_{L^2(\Omega)}\biggr) .
\end{equation}
\end{prop}
\textsc{Proof-}  
Let us first set
$$P_0=\{q\in C^2(\overline{Q_T})\ : \ q=0\hbox{ on } \Sigma_T\}.$$
The bilinear form 
$$(p,q)_P:=\iint_{Q_T}\rho^{-2} L^\star_A p\,L^\star_A q+\iint_{q_T} \rho_0^{-2}p\,q$$
where $L^\star_A q:=-q_t-\Delta q+A q$, is a scalar product on $P_0$ (see \cite{EFC-MUNCH-SEMA}). The completion $P$  of $P_0$ for the norm $\|\cdot\|_P$ associated to this scalar product is a Hilbert space and the following result proved in \cite{fursikov-imanuvilov} holds.

\begin{lemma}\label{carleman}
There exists $C=C(\Omega,\omega,T,\|A\|_\infty)>0$ such that one has the following Carleman estimate, for all $p\in P$ :
\begin{equation}\label{Carleman-ine}
\iint_{Q_T} \Big(\rho_1^{-2}|\nabla p|^2+\rho_0^{-2}|p|^2\Big)\le C\|p\|_P^2.\end{equation}
\end{lemma}
\begin{remark}
We denote by $P$ (instead of $P_A$) the completion of $P_0$ for the norm $\|\cdot\|_P$ since $P$ does not depend on $A$ (see \cite{EFC-AM}).
\end{remark}

\begin{lemma}\label{observabilite}
There exists $C=C(\Omega,\omega,T,\|A\|_\infty)>0$ such that one has the following observability inequality, for all $p\in P$ :
\begin{equation}\label{observabilite-esti}
\|p(\cdot,0)\|_{L^2(\Omega)}\le C\|p\|_P.\end{equation}
\end{lemma}
\textsc{Proof-}  
From the definition of $\rho_0$, $\rho_1$ and $\rho_2$, 
$P\xhookrightarrow{} H^1(0,\frac T2;L^2(\Omega)) \xhookrightarrow{}  C([0,\frac T2]; L^2(\Omega))$
where each imbedding is continuous. The result follows from Lemma \ref{carleman}. $\hfill\Box$

\begin{lemma}\label{solution-p}
There exists $p\in P$ unique solution of 
\begin{equation}\label{solution-p1}
(p,q)_P=\int_\Omega z_0q(0)+\int_0^T\langle \rho_2 B,\rho_2^{-1}q\rangle_{H^{-1}(\Omega)\times H^1_0(\Omega)},\quad \forall q\in P.
\end{equation}
This solution satisfies the following estimate :
$$\|p\|_P\le C \biggl(\Vert \rho_2\,B\Vert_{L^2(0,T;H^{-1}(\Omega))} + \Vert z_0\Vert_{L^2(\Omega)}\biggr) $$
where  $C=C(\Omega,\omega,T,\|A\|_\infty)>0$.
\end{lemma}
\textsc{Proof-}  The linear  map $L_1: P\to\R$, $q\mapsto \int_0^T\langle \rho_2 B,\rho_2^{-1}q\rangle_{H^{-1}(\Omega)\times H^1_0(\Omega)}$ is continuous. Indeed, for all $q\in P$
$$\Big|\int_0^T\langle \rho_2 B,\rho_2^{-1}q\rangle_{H^{-1}(\Omega)\times H^1_0(\Omega)}\Big|\le \Big(\int_0^T\|\rho_2B\|^2_{H^{-1}(\Omega)}\Big)^{1/2}
\Big(\int_0^T\|\rho_2^{-1}q\|^2_{H^{1}_0(\Omega)}\Big)^{1/2}$$
and a.e. in $(0,T)$ $\|\rho_2^{-1}q\|^2_{H^{1}_0(\Omega)}=\|\rho_2^{-1}q\|^2_{L^2(\Omega)}+\|\nabla(\rho_2^{-1}q)\|^2_{L^2(\Omega)^d}$.
But since $\rho_0\le T \rho_2$ a.e. $t$ in $(0,T)$
$$\|\rho_2^{-1}q\|^2_{L^2(\Omega)}\le \frac1{T^2}\|\rho_0^{-1}q\|^2_{L^2(\Omega)}, \,\, a.e. \ t \in (0,T).$$
Moreover 
$$\nabla(\rho_2^{-1}q)=\nabla(\rho_2^{-1})q+\rho_2^{-1}\nabla q=-\frac{s\nabla \beta(x)}{\ell(t)(T-t)^{1/2}}\rho^{-1}+\rho_2^{-1}\nabla q$$
and thus, since $\rho_1\le T^{1/2}\rho_2$ a.e. $t$ in $(0,T)$:
$$\begin{aligned}
\|\nabla(\rho_2^{-1}q)\|^2_{L^2(\Omega)^d}
&\le \biggl\|\frac{s\nabla \beta(x)}{\ell(t)(T-t)^{1/2}}\rho^{-1}q\biggl\|^2_{L^2(\Omega)^d}+\| \rho_2^{-1}\nabla q\|^2_{L^2(\Omega)^d}\\
&\le C(\Omega,\omega,T, \Vert A\Vert_{\infty})\big(\| \rho_0^{-1}q\|^2_{L^2(\Omega)}+\| \rho_1^{-1}\nabla q\|^2_{L^2(\Omega)}\big).
\end{aligned}$$
We then deduce that, a.e. in $(0,T)$
$$\|\rho_2^{-1}q\|^2_{H^{1}_0(\Omega)}\le C(\Omega, \omega ,T, \Vert A\Vert_{\infty})\big(\| \rho_0^{-1}q\|^2_{L^2(\Omega)}+\| \rho_1^{-1}\nabla q\|^2_{L^2(\Omega)^d}\big)$$
and from the Carleman estimate (\ref{Carleman-ine}) that
$$\Big(\int_0^T\|\rho_2^{-1}q\|^2_{H^{1}_0(\Omega)}\Big)^{1/2}\le C(\Omega,\omega,T, \Vert A\Vert_{\infty})\|q\|_P$$
and therefore
$$\Big|\int_0^T\langle \rho_2 B,\rho_2^{-1}q\rangle_{H^{-1}(\Omega)\times H^1_0(\Omega)}\Big|\le C(\Omega, \omega ,T, \Vert A\Vert_{\infty})\Big(\int_0^T\|\rho_2B\|^2_{H^{-1}(\Omega)}\Big)^{1/2}
\|q\|_P.$$
Thus $L_1$ is continuous.

From (\ref{observabilite-esti}) we easily deduce that the linear map $L_2:  P\to\R$, $q\mapsto \int_\Omega z_0q(0)$ is continuous.
Using Riesz's theorem, we conclude that there exists exactly one solution $p\in P$ of  (\ref{solution-p1}). $\hfill\Box$

\

Let us now introduce the convex set
$$
C(z_0,T)=\biggl\{(z,v) :  \rho z\in L^2(Q_T),\  \rho_0v\in L^2(q_T),\\
\ (z,v) \hbox{ solves } (\ref{heat_z})-(\ref{heat_z1})  \hbox{ in the transposition sense}\biggr\}
$$
that is  $(z,v)$ is solution of
$$
\iint_{Q_T}zL^\star_A q=\iint_{q_T} vq+\int_\Omega z_0q(0)+\int_0^T\langle  B,q\rangle_{H^{-1}(\Omega)\times H^1_0(\Omega)},\quad \forall q\in P.
$$
Let us remark that if $(z,v)\in C(z_0,T)$, then since $z_0\in L^2(\Omega)$, $v\in L^2(q_T)$ and $B\in L^2(0,T;H^{-1}(\Omega))$, $z$ must coincide with the unique weak solution of  (\ref{heat_z}) associated to $v$.

We can now claim that $C(z_0,T)$ is a non empty. Indeed we have :


\begin{lemma}\label{solution-controle}
Let  $p\in P$ defined in Lemma \ref{solution-p} and $(z,v)$  defined by
\begin{equation}\label{zv}
z=\rho^{-2}L^\star_A p\quad \hbox{ and }\quad v=-\rho_0^{-2} p|_{q_T}.\end{equation}
Then $(z,v)\in C(z_0,T)$ and satisfies the following estimate
\begin{equation}\label{estimation-z-v1}
\Vert \rho\, z\Vert_{L^2(Q_T)}+ \Vert \rho_0\,v\Vert_{L^2(q_T)} \leq C \biggl(\Vert \rho_2\,B\Vert_{L^2(0,T;H^{-1}(\Omega))} + \Vert z_0\Vert_{L^2(\Omega)}\biggr) 
\end{equation}
where  $C=C(\Omega,\omega,T,\|A\|_\infty)>0$.
\end{lemma}
\textsc{Proof-}  Let us prove that $(z,v)$ belongs to $C(z_0,T)$. From the definition of $P$,  $\rho  z\in L^2(Q_T)$ and $\rho _0v\in L^2(q_T)$ and from the definition of $\rho$, $\rho_0$, $\rho_2$, $z\in L^2(Q_T)$ and $v\in L^2(q_T)$. In  view of (\ref{solution-p1}), $(z,v)$ is solution of
\begin{equation}
\iint_{Q_T}zL^\star_A q=\iint_{q_T} vq+\int_\Omega z_0q(0)+\int_0^T\langle \rho_2 B,\rho_2^{-1}q\rangle_{H^{-1}(\Omega)\times H^1_0(\Omega)},\quad \forall q\in P
\end{equation}
that is, since from the definition of $\rho_2$, $B\in L^2(0,T;H^{-1}(\Omega))$ and  $\int_0^T\langle \rho_2 B,\rho_2^{-1}q\rangle_{H^{-1}(\Omega)\times H^1_0(\Omega)}=\int_0^T\langle  B, q\rangle_{H^{-1}(\Omega)\times H^1_0(\Omega)}$, $z$ is the solution of  (\ref{heat_z}) associated to $v$ in the transposition sense. Thus $C(z_0,T)\not=\emptyset$. $\hfill\Box$

Let us now consider the following extremal problem, introduced by Fursikov and Imanuvilov \cite{fursikov-imanuvilov}
\begin{equation}
\left\{\begin{aligned}\label{extrema-problem}
&\hbox{ Minimize } J(z,v)=\frac12\|(z,v)\|^2_{L^2(\rho^2;Q_T)\times L^2(\rho_0^2; q_T)}=\frac12\iint_{Q_T}\rho^2|z|^2+\frac12\iint_{q_T}\rho_0^2|v|^2\\
&\hbox{ Subject to } (z,v)\in C(z_0,T).
\end{aligned}\right.
\end{equation}
Then $(z,v)\mapsto J(z,v)$ is clearly strictly convex and continuous on $L^2(\rho^2;Q_T)\times L^2(\rho_0^2; q_T)$. Therefore  (\ref{extrema-problem}) possesses at most a unique solution in $C(z_0,T)$. More precisely we have :

\begin{prop}\label{minimal-control}
 $(z,v)\in C(z_0,T)$ defined in Lemma \ref{solution-controle} is the unique solution of (\ref{extrema-problem}).

\end{prop}
\textsc{Proof-}  Let $(y,w)\in C(z_0,T)$. Since $J$ is convex and differentiable on $L^2(\rho^2;Q_T)\times L^2(\rho_0^2; q_T)$ we have :
$$\begin{aligned}
J(y,w)
&\ge J(z,v)+\iint_{Q_T} \rho^2 z(y-z)+\iint_{q_T}\rho_0^2v(w-v)\\
&=J(z,v)+\iint_{Q_T} L^\star p (y-z)-\iint_{q_T}p(w-v)\\
&=J(z,v)
\end{aligned}$$
$y$  being  the solution of  (\ref{heat_z}) associated to $w$ in the transposition sense.
Hence $(z,v)$ solves (\ref{extrema-problem}).

To finish the proof of Proposition \ref{controllability_result}, it suffices to prove that $(z,v)$ satisfies the estimate (\ref{estimation2}).
Since $z$ is a weak solution of (\ref{heat_z}) associated to $v$, $z\in L^2(0,T;H^1_0(\Omega))$ and $z_t\in L^2(0,T;H^{-1}(\Omega))$.
Multiplying (\ref{heat_z}) by $\rho_1^2z$ and integrating by part we obtain, a.e. $t$ in $(0,T)$

$$
\begin{aligned}
\frac12\partial_t \int_\Omega |z|^2\rho_1^2&-\int_\Omega |z|^2\rho_1 \partial_t \rho_1+\int_\Omega \rho_1^2|\nabla z|^2+2\int_\Omega \rho_1z\nabla \rho_1\cdot\nabla z +\int_\Omega \rho_1^2Azz\\
&=\int_\omega v\rho_1^2 z+\langle B,\rho_1^2 z\rangle_{H^{-1}(\Omega)\times H^1_0(\Omega)}.
\end{aligned}
$$
But $\partial_t\rho_1=-\rho-(T-t)\frac{s\beta \ell'(t)}{\ell(t)^2}\rho$, so that
$$|\int_\Omega |z|^2\rho_1 \partial_t \rho_1\Big|\le C(\Omega,\omega,T, \|A\|_\infty)\int_\Omega |\rho z|^2.$$
Since $\nabla\rho_1=(T-t)\nabla\rho= (T-t)\frac{s\nabla\beta  }{\ell(t)}\rho$ we have
$$\Big|\int_\Omega \rho_1z\nabla \rho_1\cdot\nabla z \Big|\le  C(\Omega, \omega ,T, \|A\|_\infty)\Big(\int_\Omega |\rho_1\nabla z|^2\Big)^{1/2}\Big(\int_\Omega |\rho z|^2\Big)^{1/2}.$$
The following estimates also hold
$$\Big|\int_\Omega \rho_1^2 Azz\Big|\le C(T, \|A\|_\infty)\int_\Omega |\rho z|^2,$$
$$\Big|\int_\omega v\rho_1^2z\Big|\le T^{1/2}\Big|\int_\omega \rho_0 v\rho z\Big|\le T^{1/2}\Big(\int_\omega|\rho_0 v|^2\Big)^{1/.2}\Big(\int_\Omega |\rho z|^2\Big)^{1/2}$$
and
$$\begin{aligned}
|\langle B,\rho_1^2 z\rangle_{H^{-1}(\Omega)\times H^1_0(\Omega)}|
&=|\langle \rho_1B,\rho_1 z\rangle_{H^{-1}(\Omega)\times H^1_0(\Omega)}|
\le \|\rho_1 B\|_{H^{-1}(\Omega)}\|\rho_1 z\|_{H^1_0(\Omega)}\\
&\le C(\Omega, \omega, T, \|A\|_\infty) \|\rho_2 B\|_{H^{-1}(\Omega)}\big(\|\rho  z\|_{L^2(\Omega)}+\|\rho_1 \nabla z\|_{L^2(\Omega)^d}\big).
\end{aligned}
$$
Thus we easily obtain that 
$$\partial_t \int_\Omega\rho_1^2 |z|^2+\int_\Omega \rho_1^2|\nabla z|^2\le C(\Omega, \omega ,T, \|A\|_\infty)  \biggl( \|\rho_2 B\|^2_{H^{-1}(\Omega)}+\int_\Omega\rho^2 |z|^2+\int_\omega|\rho_0 v|^2\biggr) $$
and therefore, using (\ref{estimation-z-v1}), for all $t\in[0,T]$ : 
$$\Big( \int_\Omega\rho_1^2 |z|^2\Big)(t)+\iint_{Q_t} \rho_1^2|\nabla z|^2\le C(\Omega, \omega ,T, \|A\|_\infty)   \biggl(\Vert \rho_2\,B\Vert^2_{L^2(0,T;H^{-1}(\Omega))} + \Vert z_0\Vert^2_{L^2(\Omega)}\biggr)$$
which gives  (\ref{estimation2}) and concludes the proof of Proposition \ref{controllability_result}.

\section{The least-squares method and its analysis}\label{sec_LS}

For any $s\in [0,1]$, we define the space 
$$
W_s=\biggl\{g\in {\mathcal{C}}(\R), \ g(0)=0, \ g^\prime\in L^\infty(\mathbb{R}), \sup_{a,b\in \mathbb{R}, a\neq b} \frac{\vert g^\prime(a)-g^\prime(b)\vert}{\vert a-b\vert^s}< \infty\biggr\}.
$$ 
The case $s=0$ reduces to $W_0=\{g\in {\mathcal{C}}(\R),\ g(0)=0, \ g^\prime\in L^\infty(\mathbb{R})\}$ while the case $s=1$ corresponds to $
W_1=\{g\in {\mathcal{C}}(\R),\ g(0)=0, \ g^\prime\in L^\infty(\mathbb{R}), g^{\prime\prime}\in L^\infty(\mathbb{R})\}$.

In the sequel, we shall assume that there exists one $s\in (0,1]$ for which the nonlinear function $g$ belongs to $W_s$.  Remark that $g\in W_s$ for some $s\in [0,1]$
satisfies hypotheses (\ref{cond_f_5}) and (\ref{asymptotic_g}). We shall also assume that $u_0\in L^2(\Omega)$.

\subsection{The least-squares method}

We introduce the vectorial space $\mathcal{A}_0$ as follows 
\begin{equation}
\begin{aligned}
\mathcal{A}_0=
\biggl\{(y,f): \ &\rho\, y\in L^2(Q_T), \  \rho_1\, \nabla y\in L^2(Q_T)^d, \  \rho_0 f\in L^2(q_T), \\
& \rho_2(y_t - \Delta y-f\,1_{\omega} )\in L^2(0,T;H^{-1}(\Omega)),\  y(\cdot,0)=0\ \textrm{in}\ \Omega,\ y=0 \ \textrm{on}\  \Sigma_T\biggr\}
\end{aligned}
\end{equation}
where $\rho$, $\rho_2$, $\rho_1$  and $\rho_0$ are defined in \eqref{def_weight}.
Since $L^2(0,T;H^{-1}(\Omega))$ is also a Hilbert space, $\mathcal{A}_0$ endowed with the following scalar product
$$
\begin{aligned}
\big((y,f),(\overline{y},\overline{f})\big)_{\mathcal{A}_0}=\big(\rho y,\rho \overline{y}\big)_2
&+\big(\rho_1 \nabla y,\rho_1 \nabla \overline{y}\big)_2+ \big(\rho_0f,\rho_0 \overline{f}\big)_2\\
&+\big(\rho_2 (y_t-\Delta y-f\,1_{\omega}),\rho_2 (\overline{y}_t-\Delta \overline{y}-\overline{f}\,1_{\omega})\big)_{L^2(0,T;H^{-1}(\Omega))}
\end{aligned}
$$
is a Hilbert space. The corresponding norm is $\Vert (y,f)\Vert_{\mathcal{A}_0}=\sqrt{((y,f),(y,f))_{\mathcal{A}_0}}$. We also consider the convex space 
\begin{equation}
\begin{aligned}
\mathcal{A}=
\biggl\{(y,f):\ & \rho\, y\in L^2(Q_T), \  \rho_1\, \nabla y\in L^2(Q_T)^d, \  \rho_0 f\in L^2(q_T), \\
& \rho_2(y_t - \Delta y-f\,1_{\omega} )\in L^2(0,T;H^{-1}(\Omega)),\  y(\cdot,0)=u_0\ \textrm{in}\ \Omega,\ y=0 \ \textrm{on}\ \Sigma_T\biggr\}
\end{aligned}
\end{equation}
so that we can write $\mathcal{A}=(\overline{y},\overline{f})+\mathcal{A}_0$ for any element $(\overline{y},\overline{f})\in \mathcal{A}$. We endow $\mathcal{A}$ with the same norm. Clearly, if $(y,f)\in \mathcal{A}$, then  $y\in C([0,T];L^2(\Omega))$ and since $\rho\, y\in L^2(Q_T)$, then $y(\cdot,T)=0$. The null controllability requirement is therefore incorporated in the spaces $\mathcal{A}_0$ and $\mathcal{A}$. 

For any fixed $(\overline{y},\overline{f})\in \mathcal{A}$, we can now consider the following extremal problem :
\begin{equation}
\label{extremal_problem}
\min_{(y,f)\in \mathcal{A}_0} E(\overline{y}+y,f+\overline{f})
\end{equation}
where $E:\mathcal{A}\to \mathbb{R}$ is defined as follows 
\begin{equation}
E(y,f):=\frac{1}{2}\biggl\Vert \rho_2\biggl(y_t-\Delta y + g(y)-f\,1_{\omega} \biggr)\biggr\Vert^2_{L^2(0,T;H^{-1}(\Omega))}
\end{equation}
justifying the least-squares terminology we have used. 

Let us remark that, if $g\in W_s$ for one $s\geq 0$, then $g$ is Lipschitz and thus, since $g(0)=0$,  there exists $K>0$ such that $|g(\xi)|\le K|\xi|$ for all $\xi\in\R$. Consequently, $\rho_2g(y)\in L^2(Q_T)$ (and then $\rho_2g(y)\in L^2(0,T;H^{-1}(\Omega))$) since
$$\|\rho_2g(y)\|_{L^2(Q_T)}=\|(\rho_2\rho^{-1})\rho g(y)\|_{L^2(Q_T)}=\|(T-t)^{1/2}\rho g(y)\|_{L^2(Q_T)}\le T^{1/2}K\|\rho y\|_{L^2(Q_T)}.$$

Since any $g\in W_s$ satisfies hypotheses (\ref{cond_f_5}) and (\ref{asymptotic_g}), the controllability result of Theorem \ref{nullcontrollheatplus} given in \cite{EFC-EZ} implies the existence of at least one pair $(y,f)\in \mathcal{A}$ such that $E(y,f)=0$. The extremal problem (\ref{extremal_problem})
admits therefore solutions. Conversely, any pair $(y,f)\in \mathcal{A}$ for which $E(y,f)$ vanishes is a controlled pair of (\ref{heat-NL}). In this sense, the functional $E$ is a so-called error functional which measures the deviation of $(y,f)$ from being a solution of the underlying nonlinear equation. 
We emphasize that the $L^2(0,T;H^{-1}(\Omega))$ norm in $E$ indicates that we are looking for weak solutions of the parabolic equation (\ref{heat-NL}). We refer to \cite{lemoinemunch_time} where a similar so-called weak least-squares method is employed to approximate the solutions of the unsteady Navier-Stokes equation.  

A practical way of taking a functional to its minimum is through some clever use of descent  directions, i.e the use of its derivative. In doing so, the presence of local minima is always something that  may dramatically spoil the whole scheme. The unique structural property that discards this possibility is the strict convexity of the functional $E$. However, for nonlinear equation like (\ref{heat-NL}), one cannot expect this property to hold for the functional $E$. Nevertheless, we insist in that one may construct a particular minimizing sequence which cannot converge except to a global minimizer leading $E$ down to zero. 

In order to construct such minimizing sequence, we look, for any $(y,f)\in \mathcal{A}$, for a pair $(Y^1,F^1)\in \mathcal{A}_0$ solution of the following formulation 
\begin{equation}
\label{heat-Y1}
\left\{
\begin{aligned}
& Y^1_t - \Delta Y^1 +  g^{\prime}(y) Y^1 = F^1 1_{\omega}+\big(y_t-\Delta y + g(y)-f\, 1_{\omega}\big) \quad \textrm{in}\quad Q_T,\\
& Y^1=0 \,\, \textrm{on}\,\, \Sigma_T, \quad Y^1(\cdot,0)=0 \,\, \textrm{in}\,\, \Omega.
\end{aligned}
\right.
\end{equation}
Since $(Y^1,F^1)\in \mathcal{A}_0$, $F^1$ is a null control for $Y^1$. We have the following property. 

\begin{prop}  Let any $(y,f)\in \mathcal{A}$. There exists a pair $(Y^1,F^1)\in \mathcal{A}_0$ solution of (\ref{heat-Y1}) which satisfies the following estimate: 
\begin{equation} \label{estimateF1Y1}
\Vert ( Y^1,F^1)\Vert_{\mathcal{A}_0} \leq C \sqrt{E(y,f)}
\end{equation}
for some $C=C(\Omega,\omega,T, \Vert g^\prime\Vert_{\infty})>0$.\end{prop}

\textsc{Proof-}  For all $(y,f)\in \mathcal{A}$ we have
 $\rho_2(y_t-\Delta y+ g(y)-f1_\omega)\in L^2(0,T;H^{-1}(\Omega))$. The existence of a null control $F^1$ is therefore given by Proposition \ref{controllability_result}. Choosing the control $F^1$ which minimizes together with the corresponding solution $Y^1$ the functional $J$ defined in Proposition \ref{controllability_result}, we get the following estimate  (since $Y^1(\cdot,0)=0$)
\begin{equation}
\label{estimateF1Y1rho}
\begin{aligned}
\Vert \rho\, Y^1\Vert_{L^2(Q_T)}+\Vert \rho_0 F^1\Vert_{L^2(q_T)} &\leq C \Vert \rho_2 (y_t-\Delta y+g(y)-f1_\omega)\Vert_{L^2(0,T;H^{-1}(\Omega))}\\
&\leq C \sqrt{E(y,f)}
\end{aligned}
\end{equation}
and
\begin{equation}
\label{estimateF1Y1rhobis}
\begin{aligned}
 \Vert \rho_1\, Y^1\Vert_{L^\infty(0,T;L^2(\Omega))}+\Vert \rho_1 \nabla Y^1\Vert_{L^2(Q_T)^d}&\leq C \Vert \rho_2 (y_t-\Delta y+g(y)-f1_\omega)\Vert_{L^2(0,T;H^{-1}(\Omega))}\\
&\leq C \sqrt{E(y,f)}
\end{aligned}
\end{equation}
for some $C=C(\Omega,\omega,T, \Vert g\Vert_{\infty})$ independent of $Y^1$, $F^1$ and $y$. Eventually, from the equation solved by $Y^1$, 
\begin{equation}
\label{boundY1}
\begin{aligned}
\Vert\rho_2(Y^1_t-\Delta Y^1-F^1\, 1_{\omega})&\Vert_{L^2(0,T;H^{-1}(\Omega))} \\
&\leq \Vert\rho_2 g^{\prime}(y)Y^1\Vert_{L^2(Q_T)}+ \Vert \rho_2(y_t-\Delta y+g(y)-f\, 1_{\omega})\Vert_{L^2(0,T;H^{-1}(\Omega))}\\
&\leq \Vert (T-t)^{1/2}g^{\prime}(y)\Vert_{\infty}\Vert \rho Y^1\Vert_{L^2(Q_T)}+ \sqrt{2E(y,f)}\\
& \leq \max\big(1,\Vert (T-t)^{1/2}g^{\prime}\Vert_{\infty}\big)C \sqrt{E(y,f)}
\end{aligned}
\end{equation}
which proves that $(Y^1,F^1)$  belongs to $\mathcal{A}_0$.
$\hfill\Box$


\begin{remark}
From (\ref{heat-Y1}), $z=y-Y^1$ is a null controlled solution satisfying 
\begin{equation}
\label{heat-Y1-bis}
\left\{
\begin{aligned}
& z_t - \Delta z +  g^{\prime}(y) z = (f-F^1) 1_{\omega} -g(y)+g^{\prime}(y)y \quad \textrm{in}\quad Q_T,\\
& z=0 \,\, \textrm{on}\,\, \Sigma_T, \quad z(\cdot,0)=u_0 \,\, \textrm{in}\,\, \Omega
\end{aligned}
\right.
\end{equation}
by the control $(f-F^1)\in L^2(\rho_0, q_T):=\{f; \rho_0 f\in L^2(q_T)\}$.
\end{remark}

\begin{remark}
We emphasize that the presence of a right hand side in (\ref{heat-Y1}), namely $y_t-\Delta y+g(y)-f\,1_{\omega}$, forces us to introduce from the beginning the weights $\rho_0$, $\rho_1$, $\rho_2$  and $\rho$ in the spaces $\mathcal{A}_0$ and 
 $\mathcal{A}$. This can be seen from the equality \eqref{solution-p1}:  since $\rho_2^{-1} q$ belongs to $L^2(0,T; H^1(\Omega))$ for all $q\in P$, we need to impose that $\rho_2 B \in L^2(0,T; H^{-1}(\Omega))$ with here $B= y_t-\Delta y+g(y)-f\,1_{\omega}$. Working with the linearized equation \eqref{NL_z} (introduced in \cite{EFC-EZ}) which does not make appear an additional right hand side, we may avoid the introduction of Carleman type weights. Actually, the authors in \eqref{NL_z} consider controls of minimal $L^\infty(q_T)$ norm. Introduction of weights allows however the characterization \eqref{solution-p1}, which is very convenient at the practical level. We refer to \cite{EFC-MUNCH-SEMA} where this is discussed at length. 
\end{remark}

The interest of the pair $(Y^1,F^1)\in \mathcal{A}_0$ lies in the following result.
\begin{prop}\label{proposition3}
Let $(y,f)\in \mathcal{A}$  and let $(Y^1,F^1)\in \mathcal{A}_0$ be a solution of (\ref{heat-Y1}). Then the derivative of $E$ at the point $(y,f)\in \mathcal{A}$ along the direction $(Y^1,F^1)$ given by $E^{\prime}(y,f)\cdot (Y^1,F^1):=\lim_{\eta\to 0,\eta\neq 0} \frac{E((y,f)+\eta (Y^1,F^1))-E(y,f)}{\eta}$ satisfies 
\begin{equation}\label{estimateEEprime}
E^{\prime}(y,f)\cdot (Y^1,F^1)=2E(y,f).
\end{equation}
\end{prop} 
\textsc{Proof-} We preliminary check that for all $(Y,F)\in\mathcal{A}_0$,  $E$  is differentiable at the point $(y,f)\in \mathcal{A}$ along the direction $(Y,F)\in \mathcal{A}_0$. For all $\lambda\in\R$, simple computations lead to the equality 

$$
\begin{aligned}
E(y+\lambda Y,f+\lambda F) =E(y,f)+ \lambda  E^{\prime}(y,f)\cdot (Y,F) + h((y,f),\lambda (Y,F))
\end{aligned}
$$
with 
\begin{equation}\label{Efirst}
E^{\prime}(y,f)\cdot (Y,F):=\biggl(\rho_2 (y_t-\Delta y+ g(y)-f\, 1_\omega),   \rho_2 (Y_t-\Delta Y+ g^\prime(y)Y-F\, 1_\omega)\biggr)_{L^2(0,T;H^{-1}(\Omega))}
\end{equation}
and
$$
\begin{aligned}
h((y,f),\lambda (Y,F))=& \lambda \biggl(\rho_2 (Y_t-\Delta Y+ g^\prime(y)Y-F\, 1_\omega),\rho_2 l(y,\lambda Y)\biggl)_{L^2(0,T;H^{-1}(\Omega))}\\
& +\frac{\lambda ^2}2\| \rho_2 (Y_t-\Delta Y+ g^\prime(y)Y-F\, 1_\omega)\|_{L^2(0,T;H^{-1}(\Omega))}^2\\
& +\biggl(\rho_2 (y_t-\Delta y+ g(y)-f\, 1_\omega),\rho_2 l(y,\lambda Y)\biggr)_{L^2(0,T;H^{-1}(\Omega))}\\
&+ \frac{1}{2}\|\rho_2 l(y,\lambda Y)\|_{L^2(0,T;H^{-1}(\Omega))}^2
\end{aligned}
$$ 
where $l(y,\lambda Y)=g(y+\lambda Y)-g(y)-\lambda g^{\prime}(y)Y$.

The application $(Y,F)\to E^{\prime}(y,f)\cdot (Y,F)$ is linear and continuous from $\mathcal{A}_0$ to $\mathbb{R}$ as it satisfies
$$
\begin{aligned}
\vert E^{\prime}&(y,f)\cdot (Y,F)\vert \\
& \leq \Vert \rho_2 (y_t-\Delta y+ g(y)-f\, 1_\omega)\Vert_{L^2(0,T;H^{-1}(\Omega))} \Vert \rho_2 (Y_t-\Delta Y+ g^\prime(y)Y-F\, 1_\omega)\Vert_{L^2(0,T;H^{-1}(\Omega))}\\
& \leq \sqrt{2E(y,f)} \biggl(\Vert \rho_2 (Y_t-\Delta Y-F\, 1_\omega)\Vert_{L^2(0,T;H^{-1}(\Omega))} + \Vert \rho_2 g^\prime(y)Y\Vert_{L^2(Q_T)} \biggr)\\
& \leq \sqrt{2E(y,f)} \biggl(\Vert \rho_2 (Y_t-\Delta Y-F\, 1_\omega)\Vert_{L^2(0,T;H^{-1}(\Omega))} + \Vert (T-t)^{1/2} g^\prime(y)\Vert_{L^{\infty}(Q_T)}\Vert\rho Y\Vert_{L^2(Q_T)} \biggr)\\
& \leq \sqrt{2E(y,f)} \max\biggl(1,\Vert (T-t)^{1/2} g^\prime\Vert_{\infty}\biggr) \Vert (Y,F)\Vert_{\mathcal{A}_0}.
\end{aligned}
$$
Similarly, for all $\lambda \in\R^\star$
$$
\begin{aligned}
|\frac{1}{\lambda }h((y,f),\lambda (Y,F))| \leq
&  \biggl(\lambda \Vert \rho_2 (Y_t-\Delta Y+ g^\prime(y)Y-F\, 1_\omega)\Vert_{L^2(0,T;H^{-1}(\Omega))} +\sqrt{2E(y,f)}\\
&\hskip 2cm +\frac12\Vert \rho_2 l(y,\lambda Y)\Vert_{L^2(0,T;H^{-1}(\Omega))}\biggl)\frac1{\lambda }\Vert \rho_2 l(y,\lambda Y)\Vert_{L^2(0,T;H^{-1}(\Omega))}\\
&+\frac{\lambda }2\Vert \rho_2 (Y_t-\Delta Y+ g^\prime(y)Y-F\, 1_\omega)\Vert_{L^2(0,T;H^{-1}(\Omega))}^2.
\end{aligned}
$$
Since $g^{\prime}\in L^\infty(\R)$ we have for a.e. $(x,t)\in Q_T$ : 
$$\rho_2|\frac{1}{\lambda }l(y,\lambda Y)|=\rho_2\Big|\frac{g(y+\lambda Y)-g(y)}{\lambda }-g'(y)Y \Big |\le 2\|g^{\prime}\|_\infty|\rho_2Y|$$
and $\rho_2Y\in L^2(Q_T)$. Moreover, for a.e. $(x,t)\in Q_T$,  $\rho_2|\frac{1}{\lambda }l(y,\lambda Y)|=\rho_2|\frac{g(y+\lambda Y)-g(y)}{\lambda }-g'(y)Y|\to 0$  as $\lambda\to 0$; 
 it follows from the Lebesgue's Theorem that
$$\frac1{\lambda }\Vert \rho_2 l(y,\lambda Y)\Vert_{L^2(Q_T)}\to 0 \hbox{ as } \lambda\to 0.$$
It is now easy to see that 
$$h((y,f),\lambda (Y,F))=o(\lambda)$$
 and that the functional $E$  is differentiable at the point $(y,f)\in \mathcal{A}$ along the direction $(Y,F)\in\mathcal{A}_0$.
Eventually, the equality (\ref{estimateEEprime}) follows from the definition of the pair $(Y^1,F^1)$ given in (\ref{heat-Y1}).
$\hfill\Box$

\vskip 0.25cm
Remark that from the equality \eqref{Efirst}, the derivative $E^{\prime}(y,f)$ is independent of $(Y,F)$. We can then define the norm $\Vert E^{\prime}(y,f)\Vert_{(\mathcal{A}_0)^{\prime}}:= \sup_{(Y,F)\in \mathcal{A}_0, (Y,F)\neq (0,0)} \frac{E^{\prime}(y,f)\cdot (Y,F)}{\Vert (Y,F)\Vert_{\mathcal{A}_0}}$ associated to $(\mathcal{A}_0)^{\prime}$, the set of the linear and continuous applications from $\mathcal{A}_0$ to $\mathbb{R}$.

\vskip 0.25cm

Combining the equality \eqref{estimateEEprime} and the inequality \eqref{estimateF1Y1}, we deduce the following estimates of $E(y,f)$ in term of the norm of $E^\prime(y,f)$.

\begin{prop}
For any  $(y,f)\in \mathcal{A}$, the inequalities holds true
$$
C_1(\Omega,\omega,T,\|g^\prime\|_\infty) \|E'(y,f)\|_{ \mathcal{A}_0'}\le \sqrt{E(y,f)}\le C_2(\Omega,\omega,T,\|g^\prime\|_\infty) \|E'(y,f)\|_{ \mathcal{A}_0'}
$$
for some constants $C_1,C_2>0$.
\end{prop} 
\textsc{Proof-} \eqref{estimateEEprime} rewrites
$E(y,f)=\frac12 E^{\prime}(y,f)\cdot (Y^1,F^1)$
where $(Y^1,F^1)\in \mathcal{A}_0$ is solution of (\ref{heat-Y1})  and therefore, with \eqref{estimateF1Y1}
$$E(y,f)\le\frac12 \|E'(y,f)\|_{ \mathcal{A}_0'} \|(Y^1,F^1)\|_{ \mathcal{A}_0}\le C(\Omega,\omega,T,\|g^\prime\|_\infty) \|E'(y,f)\|_{ \mathcal{A}_0'}\sqrt{E(y,f)}.
$$
On the other hand, for all $(Y,F)\in\mathcal{A}_0$ (see the proof of Proposition \ref{proposition3}) :
$$
\vert E^{\prime}(y,f)\cdot (Y,F)\vert \leq  \sqrt{2E(y,f)} \max\biggl(1,\Vert (T-t)^{1/2} g^\prime\Vert_{\infty}\biggr) \Vert (Y,F)\Vert_{\mathcal{A}_0}$$
and thus 
$$C_1(\Omega,\omega,T,\|g^\prime\|_\infty) \|E'(y,f)\|_{ \mathcal{A}_0'}\le \sqrt{E(y,f)}.$$
$\hfill\Box$

In particular, any \textit{critical} point $(y,f)\in \mathcal{A}$ for $E$ (i.e. for which $E^\prime(y,f)$ vanishes) is a zero for $E$, a pair solution of the controllability problem. 
In other words, any sequence $(y_k,f_k)_{k>0}$ satisfying $\Vert E^\prime(y_k,f_k)\Vert_{(\mathcal{A}_0)^\prime}\to 0$ as $k\to \infty$ 
is such that $E(y_k,f_k)\to 0$ as $k\to \infty$. We insist that this property does not imply the convexity of the functional $E$ (and \textit{a fortiori} the strict convexity of $E$, which actually does not hold here in view of the multiple zeros for $E$) but show that a minimizing sequence for $E$ can not be stuck in a local minimum. Far from the zeros of $E$, in particular, when $\Vert (y,f)\Vert_{\mathcal{A}}\to \infty$, the right hand side inequality indicates that $E$ tends to be convex. On the other side, the left inequality indicates the functional $E$ is flat around its zero set. As a consequence, gradient based minimizing sequences may achieve a very low rate of convergence (we refer to \cite{AM-PP-2014}
and also \cite{lemoine-Munch-Pedregal-AMO-20} devoted to the Navier-Stokes equation where this phenomenon is observed).

\subsection{A strongly converging minimizing sequence for $E$}

 We now examine the convergence of an appropriate sequence $(y_k,f_k)\in \mathcal{A}$. In this respect, we observe that the equality (\ref{estimateEEprime}) shows that $-(Y^1,F^1)$ given by the solution of (\ref{heat-Y1}) is a descent direction for the functional $E$. Therefore, we can define at least formally, for any $m\geq 1$, a minimizing sequence $(y_k,f_k)_{k\in\mathbb{N}}$ as follows: 
\begin{equation}
\label{algo_LS_Y}
\left\{
\begin{aligned}
&(y_0,f_0) \in \mathcal{A}, \\
&(y_{k+1},f_{k+1})=(y_k,f_k)-\lambda_k (Y^1_k,F_k^1), \quad k\ge 0, \\
& \lambda_k= argmin_{\lambda\in[0,m]} E\big((y_k,f_k)-\lambda (Y^1_k,F_k^1)\big)    
\end{aligned}
\right.
\end{equation}
where $(Y^1_k,F_k^1)\in \mathcal{A}_0$ is such that $F^1_k$ is a null control for $Y^1_k$, solution of 

\begin{equation}
\label{heat-Y1k}
\left\{
\begin{aligned}
& Y^1_{k,t} - \Delta Y^1_k +  g^{\prime}(y_k) Y^1_k = F^1_k 1_{\omega}+ (y_{k,t}-\Delta y_k+g(y_k)-f_k 1_\omega) \quad \textrm{in}\quad Q_T,\\
& Y_k^1=0 \,\, \textrm{on}\,\, \Sigma_T, \quad Y_k^1(\cdot,0)=0 \,\, \textrm{in}\,\, \Omega
\end{aligned}
\right.
\end{equation}
and minimizes the functional $J$ defined in Proposition \ref{controllability_result}.  The direction $Y^1_k$ vanishes when $E$ vanishes. 

We first perform the analysis assuming the non linear function $g$ in $W_1$, notably that $g^{\prime\prime}\in L^{\infty}(\mathbb{R})$ (the derivatives here are in the sense of distribution). We first prove the following lemma. 

\begin{lemma}\label{estimW1lemma}
Assume $g\in W_1$. Let $(y,f)\in \mathcal{A}$ and $(Y^1,F^1)\in \mathcal{A}_0$ defined by (\ref{heat-Y1}). For any $\lambda\in \mathbb{R}$ and $k\in \mathbb{N}$, the following estimate holds 
\begin{equation}
E\big((y,f)-\lambda (Y^1,F^1)\big)  \leq    E(y,f) \biggl(\vert 1-\lambda\vert  +\lambda^2 C(\Omega,\omega,T,\|g'\|_\infty)\Vert g^{\prime\prime}\Vert_{\infty} \sqrt{E(y,f)}\biggr)^2. \label{estimW1}
\end{equation}
\end{lemma} 
\textsc{Proof-} With $g\in W_1$, we write that 
\begin{equation}\label{majoration-l}
|l(y,-\lambda Y^1)|=|g(y-\lambda Y^1)-g(y)+\lambda g^{\prime}(y) Y^1|\le  \frac{\lambda^2}{2} \|g^{\prime \prime}\|_\infty(Y^1)^2
\end{equation}
and obtain that
\begin{equation}\label{E_expansionb}
\begin{aligned}
&2 E\big((y,f)-\lambda (Y^1,F^1)\big) \\
&=\biggl\Vert \rho_2\big(y_t-\Delta y_k+g(y)-f\,1_\omega\big)-\\
&\hspace{3cm}\lambda \rho_2 \big(Y^1_t-\Delta Y^1+g^\prime(y)Y^1-F\,1_\omega\big)+\rho_2l(y,-\lambda Y^1)\biggr\Vert^2_{L^2(0,T;H^{-1}(\Omega))}\\
& =\biggl\Vert \rho_2(1-\lambda)\big(y_t-\Delta y+g(y)-f\,1_\omega\big)+\rho_2l(y,-\lambda Y^1)\biggr\Vert^2_{L^2(0,T;H^{-1}(\Omega))}\\
&\le \Big(\bigl\Vert \rho_2(1-\lambda)\big(y_t-\Delta y+g(y)-f\,1_\omega\big)\bigr \Vert_{L^2(0,T;H^{-1}(\Omega))} +\bigl \Vert\rho_2l(y,-\lambda Y^1)\bigr\Vert_{L^2(0,T;H^{-1}(\Omega))}\Big)^2\\
&\leq 2 \biggl( \vert 1-\lambda\vert \sqrt{E(y,f)}+\frac{\lambda^2}{2\sqrt{2}}\Vert g^{\prime\prime}\Vert_\infty \Vert \rho_2 (Y^1)^2\Vert_{L^2(0,T;H^{-1}(\Omega))}\biggr)^2.
\end{aligned}
\end{equation}

For $d=3$ (similar estimates hold for $d=1$ and $d=2$), using the continuous embedding of $L^{6/5}(\Omega)$ into $H^{-1}(\Omega)$, we have:
$$
\begin{aligned}
\Vert \rho_2 (Y^1)^2\Vert_{L^2(0,T;H^{-1}(\Omega))}^2
&\le C(\Omega)  \Vert \rho_2 (Y^1)^2\Vert_{L^2(0,T;L^{6/5}(\Omega))}^2 \\
&\le C(\Omega) \int_0^T\Vert \rho_2 Y^1\Vert_{L^3(\Omega)}^2\Vert   Y^1\Vert_{ L^2(\Omega)}^2\\
&\le C(\Omega) \int_0^T\Vert \rho Y^1\Vert_{L^2(\Omega)}\Vert \rho_1 Y^1\Vert_{L^6(\Omega)}\Vert   Y^1\Vert_{ L^2(\Omega)}^2\\
&\le C(\Omega) \int_0^T\Vert \rho Y^1\Vert_{L^2(\Omega)}\Vert \nabla(\rho_1 Y^1)\Vert_{L^2(\Omega)^d}\Vert   Y^1\Vert_{ L^2(\Omega)}^2.
\end{aligned}
$$
From the definition of $\rho$ and $\rho_1$ we have $\nabla \rho_1=\frac{s\nabla \beta}{\ell(t)(T-t)}\rho_1=\frac{s\nabla \beta}{\ell(t)}\rho$ and therefore a.e.  $t$ in $(0,T)$
$$
\begin{aligned}
\Vert \nabla(\rho_1 Y^1)\Vert_{L^2(\Omega)^d}
&\le  \Vert \nabla(\rho_1) Y^1\Vert_{L^2(\Omega)^d}+\Vert \rho_1\nabla Y^1\Vert_{L^2(\Omega)^d}\\
&\le C(\Omega,\omega,T,\Vert g^\prime\Vert_\infty)  \Vert \rho  Y^1\Vert_{L^2(\Omega)}+\Vert \rho_1\nabla Y^1\Vert_{L^2(\Omega)^d}\end{aligned}
$$
and thus
$$
\begin{aligned}
\Vert \rho_2 (Y^1)^2\Vert_{L^2(0,T;H^{-1}(\Omega))}^2
\le C(\Omega,\omega,T,\Vert g^\prime\Vert_\infty)
& \Vert \rho_1 Y^1\Vert_{L^\infty(0,T;L^2(\Omega))}^2\Vert  \rho Y^1\Vert_{L^2(Q_T)}\\
&\times \big(\Vert  \rho Y^1\Vert_{L^2(Q_T)}+\Vert \rho_1 \nabla Y^1\Vert_{L^2(Q_T)^d} \big).
\end{aligned}
$$

Using (\ref{estimateF1Y1rho}) and (\ref{estimateF1Y1rhobis}), we obtain
\begin{equation}\label{majoration-rhoY^2}
\Vert \rho_0 (Y^1)^2\Vert_{L^2(0,T;H^{-1}(\Omega))}^2\le {C(\Omega,\omega,T,\|g'\|_\infty)}E(y,f)^2,
\end{equation}
from which we get \eqref{estimW1}.

Proceeding as in \cite{JL-AM}, we are now in position to prove the following convergence result for the sequence $(E(y_k,f_k))_{(k\geq 0)}$.
\begin{prop}\label{convergenceEk}
Assume $g\in W_1$. Let $(y_k,f_k)_{k\in\mathbb{N}}$ be the sequence defined by (\ref{algo_LS_Y}). Then $E(y_k,f_k)\to 0$ as $k\to \infty$. Moreover, there exists $k_0\in \mathbb{N}$ such that 
the sequence $(E(y_k,f_k))_{k\geq k_0}$ decays quadratically.
\end{prop}
\textsc{Proof-} 
We define the polynomial $p_k$ as follows 
$$
p_k(\lambda)= \vert 1-\lambda\vert + \lambda^2 c_1 \sqrt{E(y_k,f_k)} \quad  \hbox{ where }  \quad c_1:=  C(\Omega,\omega,T,\|g'\|_\infty) \Vert g^{\prime\prime}\Vert_{\infty}.
$$ 
Lemma \ref{estimW1lemma} with $(y,f)=(y_k,f_k)$ allows to write that 
\begin{equation}
c_1 \sqrt{E(y_{k+1},f_{k+1})} \leq c_1 \sqrt{E(y_k,f_k)}  p_k(\widetilde{\lambda_k}), \quad \forall k\geq 0 \label{decreaseEk}
\end{equation}
with $p_k(\widetilde{\lambda_k}):=\min_{\lambda\in[0,m]}p_k(\lambda)$.

If $c_1\sqrt{E(y_0,f_0)}< 1$ (and thus $c_1\sqrt{E(y_k,f_k)}<1$ for all $k\in\mathbb{N}$) then 
$$p_k(\widetilde{\lambda_k})=\min_{\lambda\in[0,m]}p_k(\lambda)\le p_k(1)=c_1\sqrt{E(y_k,f_k)}$$
and thus
\begin{equation}
\label{C1E}
c_1\sqrt{E(y_{k+1},f_{k+1})}\le \big(c_1\sqrt{E(y_k,f_k)}\big)^2
\end{equation}
implying that $c_1\sqrt{E(y_k,f_k)}\to 0$  as $k\to \infty$ with a quadratic rate.

If now $c_1\sqrt{E(y_0,f_0)}\ge 1$, we check that $I:=\{k\in \mathbb{N},\ c_1\sqrt{E(y_k,f_k)}\ge 1\}$ is a finite subset of $\mathbb{N}$. For all $k\in I$, since $c_1\sqrt{E(y_k,f_k)}\ge 1$,
$$\min_{\lambda\in[0,m]}p_k(\lambda)=\min_{\lambda\in[0,1]}p_k(\lambda)=p_k\Big(\frac{1}{2c_1\sqrt{E(y_k,f_k)}}\Big)=1-\frac{1}{4c_1\sqrt{E(y_k,f_k)}}$$
and thus, for all $k\in I$,
\begin{equation} \label{decayEunquart}
c_1\sqrt{E(y_{k+1},f_{k+1}) } \le \Big(1-\frac{1}{4c_1\sqrt{E(y_k,f_k)}}\Big)c_1\sqrt{E(y_k,f_k)}=c_1\sqrt{E(y_k,f_k) }-\frac{1}{4}.
\end{equation}
This inequality implies that the sequence $(c_1 \sqrt{E(y_k,f_k)})_{k\in \mathbb{N}}$ strictly decreases and then that the sequence $(p_k(\widetilde{\lambda_k})_{k\in\N}$ decreases as well.   Thus the sequence  $(c_1 \sqrt{E(y_k,f_k)})_{k\in \mathbb{N}}$ decreases to $0$ at least linearly and there exists 
$k_0\in\mathbb{N}$ such that for all $k\geq k_0$, $c_1\sqrt{E(y_k,f_k)}<1$, that is  $I$ is a finite subset of $\mathbb{N}$. Arguing as in the first case, it follows that $c_1\sqrt{E(y_k,f_k)}\to 0$ as $k\to \infty$. 
In both cases, remark that $p_k(\widetilde{\lambda_k})$ decreases with respect to $k$. $\hfill\Box$

\begin{remark}
Writing from \eqref{decayEunquart} that $c_1\sqrt{E(y_k,f_k)}\leq c_1\sqrt{E(y_0,f_0)}-\frac{k}{4}$ for all $k$ such that $c_1\sqrt{E(y_k,f_k)}\geq 1$, we obtain that 
$$
k_0\leq  \biggl\lfloor 4(c_1\sqrt{E(y_0,f_0)}-1)+1 \biggr\rfloor
$$
where $\lfloor x\rfloor$ denotes the integer part of $x\in \mathbb{R}^+$. 
\end{remark}

We also have the following convergence of the optimal sequence $\{\lambda_k\}_{k>0}$. 

\begin{lemma}\label{lambda-k-go-1}
The sequence $\{\lambda_k\}_{k>0}$ defined in (\ref{algo_LS_Y}) converges to $1$ as $k\to \infty$.
\end{lemma}
\textsc{Proof-}  In view of \eqref{E_expansionb}, we have, as long as  $E(y_k,f_k)>0$, since $\lambda_k\in[0,m]$ 
$$\begin{aligned}
(1-\lambda_k)^2
&=\frac{E(y_{k+1},f_{k+1})}{E(y_{k},f_k)}-2(1-\lambda_k)\frac{\langle \rho_2\big(y_{k,t}+\Delta y_k+g(y_k)-f_k\,1_\omega\big),\rho_2l(y_k,\lambda_k Y_k^1) \rangle_{L^2(0,T;H^{-1}(\Omega))}}{E(y_{k},f_k)}\\
& \hspace{3cm}- \frac{\bigl \Vert\rho_2l(y_k,\lambda_k Y_k^1)\bigr\Vert^2_{L^2(0,T;H^{-1}(\Omega))}}{2E(y_{k})}\\
&\le \frac{E(y_{k+1},f_{k+1})}{E(y_{k},f_k)}-2(1-\lambda_k)\frac{\langle \rho_2\big(y_{k,t}+\Delta y_k+g(y_k)-f_k\,1_\omega\big),\rho_2l(y_k,\lambda_k Y_k^1) \rangle_{L^2(0,T;H^{-1}(\Omega))}}{E(y_{k},f_k)}\\
&\le \frac{E(y_{k+1},f_{k+1})}{E(y_{k},f_k)}+2\sqrt{2}m\frac{\sqrt{E(y_k,f_k)}\|\rho_2l(y_k,\lambda_k Y_k^1) \|_{L^2(0,T;H^{-1}(\Omega))}}{E(y_{k},f_k)}\\
&\le \frac{E(y_{k+1},f_{k+1})}{E(y_{k},f_k)}+2\sqrt{2}m\frac{\|\rho_2l(y_k,\lambda_k Y_k^1) \|_{L^2(0,T;H^{-1}(\Omega))}}{\sqrt{E(y_k,f_k)}}
\end{aligned}
$$

But, from \eqref{majoration-l} and  \eqref{majoration-rhoY^2} 
$$\begin{aligned}
\|\rho_2l(y_k,\lambda_k Y_k^1) \|_{L^2(0,T;H^{-1}(\Omega))}
&\le \frac{\lambda_k^2}{2\sqrt{2}}\Vert g^{\prime\prime}\Vert_\infty \Vert \rho_2 (Y^1_k)^2\Vert_{L^2(0,T;H^{-1}(\Omega))}\\
&\le m^2\Vert g^{\prime\prime}\Vert_\infty{C(T,\Omega,\omega,\|g'\|_\infty)}E(y_k,f_k)
\end{aligned}$$
and thus 
$$
(1-\lambda_k)^2\le \frac{E(y_{k+1},f_{k+1})}{E(y_{k},f_k)}+m^2\Vert g^{\prime\prime}\Vert_\infty{C(\Omega,\omega,T,\|g'\|_\infty)}\sqrt{E(y_k,f_k)}.
$$
Consequently, since $E(y_{k},f_k)\to 0$ and $\frac{E(y_{k+1},f_{k+1})}{E(y_{k},f_k)}\to 0$, we deduce that $(1-\lambda_k)^2\to 0$. $\hfill\Box$

 We are now in position to prove the following convergence result. 
 
 \begin{theorem}\label{ths1}
 Assume $g\in W_1$. Let $(y_k,f_k)_{k\in\mathbb{N}}$ be the sequence defined by (\ref{algo_LS_Y}). Then, $(y_k,f_k)_{k\in \mathbb{N}}\to (y,f)$ in $\mathcal{A}$ where $f$ is a null control for $y$ solution of  \eqref{heat-NL}. Moreover, the convergence is quadratic after a finite number of iterates. 
\end{theorem}
\textsc{Proof-} For all $k\in \mathbb{N}$, let $F_k=- \sum_{n=0}^k \lambda_n F_n^1$ and $Y_k= \sum_{n=0}^k \lambda_n Y_n^1$. Let us prove that  $\big ((Y_k,F_k )\big)_{k\in\mathbb{N}}$ converge in $\mathcal{A}_0$, i.e. that the series $\sum \lambda_n (F_n^1,Y_n^1)$ converges in $ \mathcal{A}_0$. 
Using that $\Vert(Y_k^1, F_k^1)\Vert_{\mathcal{A}_0}\leq C \sqrt{E(y_k,f_k)}$ for all $k\in\N$ (see \eqref{estimateF1Y1}), we write
$$
\sum_{n=0}^k\lambda_n\Vert(Y_n^1, F_n^1)\Vert_{\mathcal{A}_0}\le m \sum_{n=0}^k\Vert(Y_n^1, F_n^1)\Vert_{\mathcal{A}_0}
\le C \sum_{n=0}^k\sqrt{E(y_n,f_n)}.
$$
But $\big(\sqrt{E(y_n,f_n)} \big )_{k\in\mathbb{N}}$ and $\big(p_k(\widetilde{\lambda_k})\big )_{k\in\mathbb{N}}$  are decreasing sequences so that
$$\sqrt{E(y_n,f_n)}\leq p_n(\widetilde{\lambda_n}) \sqrt{E(y_{n-1},f_{n-1})}\leq p_0(\widetilde{\lambda_0}) \sqrt{E(y_{n-1},f_{n-1})}\leq p_0(\widetilde{\lambda_0})^n \sqrt{E(y_0,f_0)}$$ so that, since $p_0(\widetilde{\lambda_0})<1$ :
$$
\sum_{n=0}^k \sqrt{E(y_n,f_n)}\leq \sqrt{E(y_0,f_0)} \frac{1-p(\widetilde{\lambda_0})^{k+1}}{1-p(\widetilde{\lambda_0})}\leq \sqrt{E(y_0,f_0)} \frac{1}{1-p(\widetilde{\lambda_0})}.
$$
We deduce that  the series $\sum_n \lambda_n (Y_n^1,F_n^1)$ is normally convergent and so convergent.
Consequently, there exists $(Y,F)\in\mathcal{A}_0$ such that $(Y_k,F_k)_{k\in\N}$ converges to $(Y,F)$ in $\mathcal{A}_0$.

Denoting $y=y_0+Y$ and $f=f_0+F$, we then have that $(y_k,f_k)_{k\in\N}=(y_0+Y_k,f_0+F_k)_{k\in\N}$ converges to $(y,f)$ in $\mathcal{A}$.

It suffices now to verify that the limit $(y,f)$ satisfies $E(y,f)=0$.
We write that $(Y^1_k, F^1_k)\in \mathcal{A}_0$ and $(y_k,f_k)\in \mathcal{A}$ solve the  
\begin{equation}
\label{heat-Y1kbis}
\left\{
\begin{aligned}
& Y^1_{k,t} - \Delta Y^1_k +  g^{\prime}(y_k)\cdot Y^1_k = F^1_k 1_{\omega}- (y_{k,t}-\Delta y_k+g(y_k)-f_k 1_\omega) \quad \textrm{in}\quad Q_T,\\
& Y_k^1=0 \,\, \textrm{on}\,\, \Sigma_T, \quad Y_k^1(\cdot,0)=0 \,\, \textrm{in}\,\, \Omega.
\end{aligned}
\right.
\end{equation}
Using that $(Y^1_k, F^1_k)$ goes to zero in $\mathcal{A}_0$ as $k\to \infty$, we pass to the limit in (\ref{heat-Y1kbis})
and get, since $g\in W_1$, that $(y,f)\in \mathcal{A}$ solves (\refeq{heat-NL}), that is $E(y,f)=0$. $\hfill\Box$

In particular, along the sequence $(y_k,f_k)_k$ defined by (\ref{algo_LS_Y}), we have the following coercivity property for $E$, which confirms the strong convergence of the sequence $(y_k,f_k)_{k>0}$. In view of the non uniqueness of the zeros of $E$, remark that this property is not true in general for all $(y,f)$ in $\mathcal{A}$.

\begin{prop}
Let $(y_k,f_k)_{k>0}$ defined by (\ref{algo_LS_Y}) and $(\overline{y},\overline{f})$ its limit. Then, there exists a positive constant $C$ such that  
\begin{equation}\label{estim_coercivity}
\Vert (\overline{y},\overline{f})-(y_k,f_k)\Vert_{\mathcal{A}_0} \leq C\sqrt{E(y_k,f_k)}, \quad \forall k>0.
\end{equation}
\end{prop}
\textsc{Proof-} We write that 
\begin{equation}
\nonumber
\begin{aligned}
\Vert (\overline{y},\overline{f})-(y_k,f_k)\Vert_{\mathcal{A}_0} & =\Vert \sum_{p=k+1}^{\infty} \lambda_p (Y^1_p,F^1_p)\Vert_{\mathcal{A}}\leq m\sum_{p=k+1}^{\infty}  \Vert (Y^1_p,F^1_p \Vert_{\mathcal{A}_0}\\ 
& \leq m\, C\sum_{p=k+1}^{\infty}  \sqrt{E(y_p,f_p)}\\
& \leq m\, C\sum_{p=k+1}^{\infty}  p_{0}(\widetilde{\lambda}_0)^{p-k}\sqrt{E(y_k,f_k)}\\
& \leq m\, C\frac{p_{0}(\widetilde{\lambda}_0)}{1-p_{0}(\widetilde{\lambda}_0)}\sqrt{E(y_{k},f_{k})}.
\end{aligned}
\end{equation}
$\hfill\Box$

We emphasize, in view of the non uniqueness of the zeros of $E$, that an estimate (similar to \eqref{estim_coercivity}) of the form  $\Vert (\overline{y},\overline{f})-(y,f)\Vert_{\mathcal{A}_0} \leq C \sqrt{E(y,f)}$ does not hold for all $(y,f)\in\mathcal{A}$. We also mention the fact that the sequence $(y_k,f_k)_{k>0}$ and its limits $(\overline{y},\overline{f})$ are uniquely determined from the initial guess $(y_0,f_0)$ and from our criterion of selection of the control $F^1$. In other words, the solution $(\overline{y},\overline{f})$ is unique up to the element $(y_0,f_0)$ and the functional $J$.

\subsection{The case $g\in W_s$, $0\leq s<1$ and additional remarks}

The results of the previous subsection devoted to the case $s=1$ still hold if we assume only that $g\in W_s$ for one $s\in (0,1)$. For any $g\in W_s$, we introduce the notation 
$\|g^\prime\|_{\widetilde{W}^{s,\infty}(\R)}:=\sup_{a,b\in \mathbb{R}, a\neq b} \frac{\vert g^\prime(a)-g^\prime(b)\vert}{\vert a-b\vert^s}$. We have the following result. 

 \begin{theorem}\label{ths}
Assume that there exists $s\in (0,1)$ such that $g\in W_s$. Let $(y_k,f_k)_{k\in\mathbb{N}}$ be the sequence defined by (\ref{algo_LS_Y}). Then, $(y_k,f_k)_{k\in \mathbb{N}}\to (y,f)$ in $\mathcal{A}$ where $f$ is a null control for $y$ solution of  \eqref{heat-NL}. Moreover, after a finite number of iterates, the rate of convergence is equal to $1+s$.
\end{theorem}
\textsc{Proof-}  We briefly sketch the proof, close to the proof of Theorem \ref{ths1} for the case $s=1$. 

-We first prove for any $(y,f)\in \mathcal{A}$ and $\lambda\in \R$ the following inequality (similar to the inequality \eqref{estimW1})
\begin{equation}
\label{estimE_WS}
E\big((y,f)-\lambda (Y^1,F^1)\big)  \leq   E(y,f)\biggl(\vert 1-\lambda\vert+ \lambda^{1+s} c_1
E(y,f)^{s/2}\biggr)^2
\end{equation}
with $c_1=C(T,\Omega,\omega,\|g'\|_\infty) \|g^\prime\|_{\widetilde{W}^{s,\infty}(\R)}$ and $(Y^1,F^1)\in \mathcal{A}_0$ the solution of (\ref{heat-Y1k}) which minimizes $J$.
For any $(x,y)\in \R^2$ and $\lambda\in \R$, we write $g(x+\lambda y)-g(x)=\int_0^\lambda y g'(x+\xi y)d\xi$ leading to 
$$
\begin{aligned}
|g(x+\lambda y)-g(x)-\lambda g'(x)y|
&\le \int_0^\lambda |y| |g'(x+\xi y)-g'(x)|d\xi\\
&\le \int_0^\lambda |y|^{1+s}|\xi|^s\frac{|g'(x+\xi y)-g'(x)|}{|\xi y|^s}d\xi\\
&\le \|g^\prime\|_{\widetilde{W}^{s,\infty}(\R)}|y|^{1+s}\frac{\lambda^{1+s}}{1+s}.
\end{aligned}
$$
It follows that 
$$|l(y,-\lambda Y^1)|=|g(y-\lambda Y^1)-g(y)+\lambda g^{\prime}(y) Y^1|\le  \|g^\prime\|_{\widetilde{W}^{s,\infty}(\R)}\frac{\lambda^{1+s}}{1+s}|Y^1|^{1+s}
$$
and
$$\begin{aligned}
\bigl \Vert\rho_2l(y,\lambda Y^1)\bigr\Vert_{L^2(0,T;H^{-1}(\Omega))}
&\le \bigl \Vert\rho_2l(y,\lambda Y^1)\bigr\Vert_{L^2(0,T;L^{6/5}(\Omega))}\\
&\le  \|g^\prime\|_{\widetilde{W}^{s,\infty}(\R)}\frac{\lambda^{1+s}}{1+s}\bigl \Vert\rho_2|Y^1|^{1+s}\bigr\Vert_{L^2(0,T;L^{6/5}(\Omega))}.
\end{aligned}
 $$
But
$$\begin{aligned}
 \bigl \Vert\rho_2|Y^1|^{1+s}\bigr\Vert_{L^2(0,T;L^{6/5}(\Omega))}^2
 &=\int_0^T \bigl \Vert\rho_2|Y^1|^{1+s}\bigr\Vert_{L^{6/5}(\Omega)}^2
 \le \int_0^T  \bigl \Vert\rho_2Y^1\bigr\Vert_{L^{3}(\Omega)}^2 \bigl \Vert |Y^1|^{s}\bigr\Vert_{L^{2}(\Omega)}^2\\
 &\le \int_0^T  \bigl \Vert\rho Y^1\bigr\Vert_{L^{2}(\Omega)} \bigl \Vert\rho_1Y^1\bigr\Vert_{L^{6}(\Omega)} \bigl \Vert Y^1\bigr\Vert_{L^{2s}(\Omega)}^{2s}\\
  &\le C(\Omega) \int_0^T  \bigl \Vert\rho Y^1\bigr\Vert_{L^{2}(\Omega)} \bigl \Vert\nabla(\rho_1Y^1)\bigr\Vert_{L^{2}(\Omega)^d} \bigl \Vert Y^1\bigr\Vert_{L^{2s}(\Omega)}^{2s}\\
  &\le C(\Omega) \bigl \Vert\rho Y^1\bigr\Vert_{L^{2}(Q_T)} \bigl \Vert\nabla(\rho_1Y^1)\bigr\Vert_{L^{2}(Q_T)^d} \bigl \Vert Y^1\bigr\Vert_{L^\infty(0,T; L^{2s}(\Omega))}^{2s}\\
&\le C(\Omega) \bigl \Vert\rho Y^1\bigr\Vert_{L^{2}(Q_T)} \bigl \Vert\nabla(\rho_1Y^1)\bigr\Vert_{L^{2}(Q_T)^d} \bigl \Vert Y^1\bigr\Vert_{L^\infty(0,T; L^{2}(\Omega))}^{2s}.
 \end{aligned}  
 $$
 Since $\Vert \nabla(\rho_1 Y^1)\Vert_{L^2(\Omega)^d} \le C(\Omega,\omega,T,\Vert g^\prime\Vert_\infty)  \Vert \rho  Y^1\Vert_{L^2(\Omega)}+\Vert \rho_1\nabla Y^1\Vert_{L^2(\Omega)^d}$, we finally get
$$
\begin{aligned}
 \bigl \Vert\rho_2|Y^1|^{1+s}\bigr\Vert_{L^2(0,T;L^{6/5}(\Omega))}^2
 &\le   C(\Omega,\omega,T,\|g'\|_\infty)
\Vert \rho Y^1\Vert_{L^{2}(Q_T)}\\
&\hskip 2cm \times\big(\Vert \rho Y^1\Vert_{L^{2}(Q_T)}+ \Vert \rho_1 \nabla Y^1\Vert_{L^{2}(Q_T)^d}\big)\Vert\rho_1 Y\Vert_{L^\infty(0,T;L^2(\Omega))}^{2s}.
\end{aligned}
$$
The first inequality of (\ref{E_expansionb}) then leads to (\ref{estimE_WS}).

- We then check that the sequence $(E(y_k,f_k))_{k\in \mathbb{N}}$ goes to zero as $k\to \infty$. We define $p_k$ as follows 
$$
p_k(\lambda)= \vert 1-\lambda\vert + \lambda^{1+s} c_1  E(y_k,f_k)^{s/2} 
$$ 
so that
\begin{equation}
 \sqrt{E(y_{k+1},f_{k+1})} \leq \sqrt{E(y_k,f_k)}  p_k(\widetilde{\lambda_k}), \quad \forall k\geq 0 \nonumber 
\end{equation}
with $p_k(\widetilde{\lambda_k})=\min_{\lambda\in[0,m]}p_k(\lambda)$.
We have $p_k(\widetilde{\lambda}_k):=\min_{\lambda\in[0,m]}p_k(\lambda)\le p_k(1)=c_1 E(y_k,f_k)^{s/2}$
and thus
$$c_2\sqrt{E(y_{k+1},f_{k+1})}\le \big(c_2\sqrt{E(y_k,f_k)}\big)^{1+s}, \qquad c_2:=c_1^{1/s}.$$

If $c_2\sqrt{E(y_0,f_0)}< 1$ (and thus $c_2\sqrt{E(y_k,f_k)}<1$ for all $k\in\mathbb{N}$) then 
the above inequality implies that $c_2\sqrt{E(y_k,f_k)}\to 0$  as $k\to \infty$. If $c_2\sqrt{E(y_0,f_0)}\ge 1$ then let $I=\{k\in \mathbb{N},\ c_2\sqrt{E(y_k,f_k)}\ge 1\}$. $I$ is a finite subset of $\mathbb{N}$;  for all $k\in I$, since $c_2\sqrt{E(y_k,f_k)}\ge 1$
$$\min_{\lambda\in[0,m]}p_k(\lambda)=\min_{\lambda\in[0,1]}p_k(\lambda)=p_k\Big(\frac{1}{(1+s)^{1/s}c_2\sqrt{E(y_k,f_k)}}\Big)=1- \frac{s}{(1+s)^{\frac1s+1}}\frac{1}{c_2\sqrt{E(y_k,f_k)}}$$
and thus, for all $k\in I$,
$$c_2\sqrt{E(y_{k+1},f_{k+1}) } \le \Big(1- \frac{s}{(1+s)^{\frac1s+1}}\frac{1}{c_2\sqrt{E(y_k,f_k)}}\Big)c_2\sqrt{E(y_k,f_k)}=c_2\sqrt{E(y_k,f_k) }-\frac{s}{(1+s)^{\frac1s+1}}.$$
This inequality implies that the sequence $(c_2 \sqrt{E(y_k,f_k)})_{k\in \mathbb{N}}$ strictly decreases and then that the sequence $(p_k(\widetilde{\lambda_k}))_{k\in\N}$ decreases as well.   Thus the sequence  $(c_2 \sqrt{E(y_k,f_k)})_{k\in \mathbb{N}}$ decreases to $0$ at least linearly and there exists 
$k_0\in\mathbb{N}$ such that for all $k\geq k_0$, $c_2\sqrt{E(y_k,f_k)}<1$, that is  $I$ is a finite subset of $\mathbb{N}$. Similarly, the optimal parameter $\lambda_k$ goes to one as $k\to \infty$.

- Using that the sequence $(E(y_k,f_k))_{k\in \mathbb{N}}$ goes to zero, we conclude exactly as in the proof of Theorem \ref{ths1}. $\hfill\Box$

On the other hand, if we assume only that $g$ belongs to $W_0$, then we can not expect the convergence of the sequence $(y_k,f_k)_{k>0}$ if $\Vert g^\prime\Vert_{\infty}$ is too large.

 \begin{remark} \label{remarque4}
Assume that $g\in W_0$. Let any $(y,f)\in \mathcal{A}$ and $(Y^1,F^1)$ the solution of \eqref{heat-Y1} which minimizes $J$. The following inequality holds : 
\begin{equation}
\nonumber
E\big((y,f)-\lambda (Y^1,F^1)\big)  \leq   E(y,f)\biggl(\vert 1-\lambda\vert+ \lambda C(\Omega,\omega,T,\|g'\|_\infty) \|g'\|_\infty \biggr)^2
\end{equation}
for all $\lambda\in \mathbb{R}$ where $C(\Omega,\omega,T,\|g'\|_\infty) \ge 0$ increases with $\|g'\|_\infty$. 
Indeed, this is a consequence of the following inequality, for all $(y,f)\in\mathcal{A}$, $(Y,F)\in\mathcal{A}_0$ : 
$$
\begin{aligned}
2 E\big((y,f)-\lambda (Y^1,F^1)\big)
& \le \Big(\bigl\Vert \rho_2(1-\lambda)\big(y_t-\Delta y+g(y)-f\,1_\omega\big)\bigr \Vert_{L^2(0,T;H^{-1}(\Omega))} \\
&\hskip 6cm+\bigl \Vert\rho_2l(y,\lambda Y^1)\bigr\Vert_{L^2(0,T;H^{-1}(\Omega))}\Big)^2\\
& \leq \biggl(|1-\lambda|\sqrt{2E(y,f)}+ 2\lambda  \Vert (T-t)^{1/2} g^\prime(y)\Vert_{L^{\infty}(Q_T)}\Vert\rho Y\Vert_{L^2(Q_T)}  \biggr)^2.
\end{aligned}
$$
As a consequence,  we get that the sequence $(E(y_k,f_k))_{k\ge 0}$ decreases to $0$ if $g$ satisfies 
$$
C(\Omega,\omega,T,\|g'\|_\infty) \|g'\|_\infty<1.
$$ 
\end{remark}
$\hfill\Box$

\begin{remark}
The estimate  \eqref{estimateF1Y1} is a key point in the convergence analysis and is independent of the choice of the functional $J$ defined by $J(Y^1,F^1)=\frac{1}{2}\Vert \rho_0 F^1\Vert^2_{L^2(q_T)}+\frac{1}{2}\Vert \rho Y \Vert^2_{L^2(Q_T)}$  (see Proposition \ref{controllability_result}) in order to select a pair $(Y^1,F^1)$ in $\mathcal{A}_0$. Thus, we may consider other weighted functionals, for instance $J(Y^1,F^1)=\frac{1}{2}\Vert \rho_0 F^1\Vert^2_{L^2(q_T)}$ as discussed in \cite{DeSouza_Munch_2016}. 
\end{remark}

\begin{remark}\label{link_DampedNewtonMethod}
If we introduce $F:\mathcal{A}\to L^2(0,T;H^{-1}(\Omega))$ by $F(y,f):=\rho^{-2}(y_t-\Delta y + g(y)-f\,1_\omega)$, we get that $E(y,f)=\frac{1}{2}\Vert F(y,f)\Vert_{L^2(0,T;H^{-1}(\Omega))}^2$ and observe that, for $\lambda_k=1$, the algorithm \eqref{algo_LS_Y} coincides with the Newton algorithm associated to the mapping $F$. This explains notably the quadratic convergence of Theorem \ref{ths1} in the case $g\in W_1$ for which we have a control of $g^{\prime\prime}$ in $L^{\infty}(Q_T)$. The optimization of the parameter $\lambda_k$ allows to get a global convergence of the algorithm and leads to the so-called damped Newton method (for $F$). Under general hypothesis, global convergence for this kind of method is achieved, with a linear rate (for instance; we refer to \cite[Theorem 8.7]{deuflhard}). As far as we know, the analysis of damped type Newton methods  for partial differential equations has deserved very few attention in the literature. We mention \cite{lemoinemunch_time, saramito} in the context of fluids mechanics.   
\end{remark}

\begin{remark}Suppose to simplify that $\lambda_k$ equals one (corresponding to the standard Newton method). Then, for each $k$, the optimal pair $(Y_k^1,F_k^1)\in \mathcal{A}_0$ is such that the element $(y_{k+1},f_{k+1})$ minimizes over $\A$ the functional $(z,v)\to J(z-y_k,v-f_k)$. Instead, we may also select the pair $(Y_k^1,F_k^1)$ such that the element $(y_{k+1},f_{k+1})$ minimizes the functional $(z,v)\to J(z,v)$. This leads to the following sequence $\{y_k,f_k\}_k$ defined by 
 \begin{equation}
\label{heat_lambdakequal1}
\left\{
\begin{aligned}
& y_{k+1,t} - \Delta y_{k+1} +  g^{\prime}(y_k) y_{k+1} = f_{k+1} 1_{\omega}+g^\prime(y_k)y_k-g(y_k), & \textrm{in}\,\, Q_T,\\
& y_k=0,  & \textrm{on}\,\, \Sigma_T, \\
& (y_{k+1}(\cdot,0), y_{k+1,t}(\cdot,0))=(u_0,u_1), & \textrm{in}\,\, \Omega.
\end{aligned}
\right.
\end{equation}
This is actually the formulation used in \cite{EFC-AM}. This formulation is different and the analysis of convergence (at least in the framework of our least-squares setting) is less direct because it is necessary to have a control of the right hand side term $g^\prime(y_k)y_k-g(y_k)$. 
\end{remark}

\begin{remark}
We emphasize that the explicit construction used here allows to recover the null controllability property of \eqref{heat-NL} for nonlinearities $g$ in $W_s$ for one $s\in (0,1]$. We do not use a fixed point argument as in \cite{EFC-EZ}. On the other hand, the conditions we make on $g$ are more restrictives that in \cite{EFC-EZ}. Eventually, it is also important to remark these additional conditions on $g$ does not imply a priori a contraction property of the operator $\Lambda$ introduced in \cite{EFC-EZ} and mentioned in the introduction. Assume $g\in W_1$. If $(y_{z^i},f_{z^i})$, $i=1,2$ are a controlled pair for the system \eqref{NL_z} minimizing the functional $J$, then the following inequality holds : 
\begin{equation}
\Vert \rho_0(f_{z^1}-f_{z^2})\Vert_{L^2(q_T)}+ \Vert \rho (y_{z^1}-y_{z^2})\Vert_{L^2(Q_T)} \leq C(\Omega, \omega,T,\Vert \tilde{g}\Vert_{\infty}) \Vert g^{\prime\prime}\Vert_{\infty}\Vert u_0\Vert_{L^2(\Omega)}\Vert z^1-z^2\Vert_{L^{\infty}(Q_T)}
\end{equation}
where $C(\Omega, \omega,T,\Vert \tilde{g}\Vert_{\infty})$ is the constant appearing in \eqref{estimation1}. In order to ensure a contraction property, we need a priori to add a smallness assumption on the data $g$ and $u_0$. 
\end{remark}

\section{Numerical illustrations}\label{sec_numeric}

We illustrate in this section our results of convergence. We first provide some practical details of the algorithm \eqref{algo_LS_Y} then discussed some experiments in the one dimensional case. 

\subsection{Approximation - Algorithm} 

Each iterate of the algorithm (\ref{algo_LS_Y}) requires the determination of the null control of $F^1_k$ for $Y^1_k$ solution of 
\begin{equation}
\label{heat-Y1kh}
\left\{
\begin{aligned}
& Y^1_{k,t} - \Delta Y^1_k +  g^{\prime}(y_k) Y^1_k = F^1_k 1_{\omega}+B_k, & \textrm{in}\,\, Q_T,\\
& Y_k^1=0,  & \textrm{on}\,\, \Sigma_T, \\
& Y_k^1(\cdot,0)=0, & \textrm{in}\,\, \Omega
\end{aligned}
\right.
\end{equation}
with $B_k:= y_{k,t}-\Delta y_k+g(y_k)-f_k 1_\omega$. From Lemma \ref{solution-controle}, the pair $(F^1_k,Y^1_k)$ which minimizes the functional $J$ is given by
$$
Y^1_k=\rho^{-2}L_{g^{\prime}(y_k)}^\star p_k, \quad F_k^1=-\rho_0^{-2}p_k \, 1_{q_T}
$$
where $p_k\in P$ solves the formulation 
\begin{equation}\label{FV-P_algo}
\iint_{Q_T} \rho^{-1}L_{g^{\prime}(y_k)}^{\star} p_k \, \rho^{-1}L_{g^{\prime}(y_k)}^{\star}\overline{p} + \iint_{q_T} \rho_0^{-1}p_k \, \rho_0^{-1}\overline{p} = \int_0^T < \rho_2 B_k, \rho_2^{-1}\overline{p}>_{H^{-1}(\Omega),H_0^1(\Omega)}dt\quad \forall \overline{p}\in P.
\end{equation}
The numerical approximation of this variational formulation (of second order in time and fourth order in space) has been discussed at length in \cite{EFC-MUNCH-SEMA}. In order, first to avoid numerical instabilities (due to the presence of exponential functions in the formulation), and second to make appear explicitly the controlled solution, we introduce the new variables
$$
m_k=\rho_0^{-1}p, \qquad z_k=\rho^{-1}L_{g^{\prime}(y_k)}^\star p_k.
$$
Since $\rho_2^{-1}p\in L^2(0,T;H_0^1(\Omega))$, we obtain notably that $\rho_2^{-1}p=\rho_2^{-1}\rho_0m=(T-t)m\in L^2(0,T; H_0^1(\Omega))$.
From \eqref{FV-P_algo}, the pair $(m_k,z_k)\in \mathcal{M}\times L^2(Q_T)$ with $\mathcal{M}:=\rho_0^{-1}P$ solves
\begin{equation}
\iint_{Q_T} z_k \, \overline{z} + \iint_{q_T} m_k \, \overline{m} = \int_0^T  < \rho_2 B_k, (T-t) \overline{m}>_{H^{-1}(\Omega),H_0^1(\Omega)}    dt \quad \forall (\overline{m},\overline{z})\in \mathcal{M}\times L^2(Q_T)
\end{equation}
subject to the constraint $z_k=\rho^{-1}L_{g^{\prime}(y_k)}^{\star}(\rho_0 m_k)$. This constraint leads to the following well-posed mixed formulation :  find $(m_k,z_k,\lambda_k)\in \mathcal{M}\times L^2(Q_T)\times L^2(Q_T)$ solution of 
\begin{equation}
\label{FV_mz}
\left\{
\begin{aligned}
&\iint_{Q_T} z_k \, \overline{z} + \iint_{q_T} m_k \, \overline{m} + \iint_{Q_T} \lambda_k \biggl(\overline{z}-\rho^{-1}L_{g^{\prime}(y)}^{\star}(\rho_0\, \overline{m})\biggr)\\
&\hspace{2cm}= \int_0^T < \rho_2 B_k, (T-t)\overline{m}>_{H^{-1}(\Omega),H_0^1(\Omega)}dt,
\quad \forall (\overline{m}, \overline{z})\in \mathcal{M}\times L^2(Q_T), \\
& \iint_{Q_T} \overline{\lambda} \biggl(z_k-\rho^{-1}L_{g^{\prime}(y_k)}^{\star}(\rho_0\, m)\biggr)=0, \quad\forall \overline{\lambda}\in L^2(Q_T).
\end{aligned}
\right.
\end{equation}
The variable $\lambda_k\in L^2(Q_T)$ is a Lagrange multiplier.  Moreover, from the unique solution $(m_k,z_k)$, we get the explicit form of the controlled pair $(Y^1_k,F^1_k)$ as follows:
$$
Y_k^1=\rho^{-1}z_k, \qquad F_k^1=-\rho_0^{-1} m_k \, 1_{q_T}.
$$

The algorithm associated to the sequence $(y_k,f_k)_{k>0}$ (see \eqref{algo_LS_Y}) may be developed as follows: given $\epsilon>0$ and $m\geq 1$,
\par\noindent
\begin{enumerate}
\item We determine the controlled pair $(y_0,f_0)$ which minimizes the functional $J$ associated to the linear case (for which $g\equiv 0$ in \eqref{heat-NL}). $(y_0,f_0)$ is given by  
$$
(y_0,f_0)=(\rho^{-1}z_0, -\rho_0^{-1} m_0 \, 1_{q_T})
$$
 where $(z_0,m_0)$ solves the formulation :
\begin{equation}
\left\{
\begin{aligned}
&\iint_{Q_T} z \, \overline{z} + \iint_{q_T} m \, \overline{m} + \iint_{Q_T} \lambda \biggl(\overline{z}-\rho^{-1}L_0^{\star}(\rho_0\, \overline{m})\biggr)= \iint_\Omega \rho_0(\cdot,0)\, u_0 \,\overline{m}(\cdot,0),\\
&\hspace{7cm}\forall (\overline{m}, \overline{z})\in \mathcal{M}\times L^2(Q_T), \\
& \iint_{Q_T} \overline{\lambda} \biggl(z-\rho^{-1}L_0^{\star}(\rho_0\, m)\biggr)=0, \quad\forall \overline{\lambda}\in L^2(Q_T).
\end{aligned}
\right.
\end{equation}
In view of Proposition \ref{controllability_result}, we check that $(y_0,f_0)$ belongs to $\A$.

\item Assume now that $(\lambda_k,f_k)$ is computed for some $k\geq 0$. We then compute $c_k\in L^2(0,T;H_0^1(\Omega))$, unique solution of  
\begin{equation}
\int_{Q_T} \nabla c_k\cdot \nabla \overline{c} =\int_0^T <\rho_2(y_{k,t}-\Delta y_k+g(y_k)-f_k\,1_\omega),\overline{c}>_{H^{-1}(\Omega),H_0^1(\Omega)}, \quad \forall \overline{c}\in L^2(0,T;H_0^1(\Omega))
\end{equation}
and then $E(y_k,f_k)=\frac{1}{2}\Vert \rho_2(y_{k,t}-\Delta y_k+g(y_k)-f_k\,1_\omega)\Vert^2_{L^2(0,T;H^{-1}(\Omega))}=\frac{1}{2}\Vert \nabla c_k\Vert^2_{L^2(Q_T)}$. 

\item If $E(y_k,f_k)< \epsilon$, the approximate controlled pair is given by $(\overline{y},\overline{f})=(y_k,f_k)$ and the algorithm stops. Otherwise, we determine  
the solution $(Y^1_k,F^1_k)=(\rho^{-1}z_k,-\rho_0^{-1}m_k\,1_{q_T})$ where $(z_k,m_k)$ solves (\ref{FV_mz}).

\item Set $(y_{k+1},f_{k+1})=(y_k,f_k)-\lambda_k (Y^1_k,F^1_k)$ where $\lambda_k$ minimizes over $[0,m]$ the scalar functional $\lambda\rightarrow  E((y_k,f_k)-\lambda (Y^1_k,F^1_k))$ defined by (see \eqref{E_expansionb})
\begin{equation}
\begin{aligned}
2 E\big((y_k,f_k)-\lambda (Y_k^1,F_k^1)\big)
=\biggl\Vert \rho_2(1-\lambda)\big(y_{k,t}-\Delta y_k+g(y_k)-f_k\,1_\omega\big)+\rho_2l(y_k,-\lambda Y_k^1)\biggr\Vert^2_{L^2(0,T;H^{-1}(\Omega))}\\
\end{aligned}
\end{equation}
with $l(y_k,-\lambda Y_k^1)=g(y_k-\lambda Y_k^1)-g(y_k)+\lambda g^{\prime}(y_k)Y_k^1$. The minimization is performed using a line search method. Return to step 2.

\end{enumerate}

We use  the conformal space-time finite element method described in \cite{EFC-MUNCH-SEMA}. We consider a regular family $\mathcal{T}=\{\mathcal{T}_h; h>0\}$ of triangulation of $Q_T$ such that $\overline{Q_T}=\cup_{K\in \mathcal{T}_h}  K$. The family $\mathcal{T}$ is indexed by $h=max_{K\in \mathcal{T}_h} diam(K)$. The variable $z_k$ and $\lambda_k$ are approximated with the space $P_h=\{p_h\in C(\overline{Q_T}); p_h\vert_K\in \mathbb{P}_1(K), \forall K\in \mathcal{T}_h\}\subset L^2(Q_T)$ where $\mathbb{P}_1(K)$ denotes the space of affine functions both in $x$ and $t$. The variable $m_k$ is approximated with the space  $V_h=\{v_h\in C^1(\overline{Q_T}); v_h\vert_K\in \mathbb{P}(K), \forall K\in \mathcal{T}_h\}\subset \mathcal{M}$ where $\mathbb{P}(K)$ denotes the Hsieh-Clough-Tocher $C^1$ element (we refer to \cite{Ciarlet78} page 356). These conformal approximation leads to a strong convergent approximation of the control and the controlled solution with respect to the parameter $h$.

\subsection{Experiments}

We present some numerical experiments in the one dimensional setting with $\Omega=(0,1)$. The control is located on $\omega=(0.1,0.3)$. We consider $T=1/2$; moreover, in order to reduce the dissipation of the solution of (\ref{heat-NL}) when $g\equiv 0$, we replace the term $-\Delta y$ in (\ref{heat-NL}) by $-\nu\Delta y$ with $\nu>0$ small, here $\nu=10^{-1}$. We consider the nonlinear even function $g$ as follows 

\begin{equation}
g(s)=\left\{
\begin{aligned}
& l(s), & s\in [-a,a],\\
& -\vert s\vert^\alpha \log^{3/2}(1+\vert s\vert), & \vert s\vert\geq a   
\end{aligned}
\right.
\end{equation}
with $a,\alpha\in (0,1)$. $l$ denotes the (even) polynomial of order two such that $l(0)=0$, $l(a)= -\vert a\vert^\alpha \log^{3/2}(1+\vert a\vert)$ and $(-\vert s\vert^\alpha \log^{3/2}(1+\vert s\vert))^{\prime}(s=a)=l^\prime(a)$.
We use in the sequel the values $a=10^{-1}$ and $\alpha=0.95$. We check that $g$ belongs to $W_1$, in particular $g^{\prime\prime}\in L^{\infty}(\mathbb{R})$ in the sense of distribution. Remark as well that $g$ is sublinear. 

As for the initial condition to be controlled, we consider simply $u_0(x)=\beta \sin(\pi x)$ parametrized by $\beta>0$.

The experiments are performed with the Freefem++ package developed at the Sorbonne university  (see \cite{Hecht2012}), very well-adapted to the space-time formulation we employ. 
The algorithm is stopped when the value $E(y_k,f_k)$ is less than $\epsilon=10^{-6}$. The optimal steps $\lambda_k$ are searched in the interval $[0,1]$.

Table \ref{tab:LS_beta10}, \ref{tab:LS_beta100} and \ref{tab:LS_beta1000} collect some norms from the sequence $(y_k,f_k)_{k\geq 0}$ defined by the algorithm (\ref{algo_LS_Y}), initialized with the linear controlled solution, for $\beta=10.$,  $\beta=10^2$ and $\beta=10^3$ respectively. We use a structured mesh composed of $20\ 000$ triangles, $10\ 201$ vertices and for which $h\approx 1.11\times 10^{-2}$.  For $\beta=10$, we observe the convergence after $4$ iterates. The optimal steps $\lambda_k$ are very close to one since $\max_{k} \vert \lambda_k-1\vert< 0.05$; consequently, the algorithm (\ref{algo_LS_Y}) provides similar results than the Newton algorithm (for which $\lambda_k=1$ for all $k$). For $\beta=10^2$, the convergence remains fast and is reached after 8 iterates. We can observe that some optimal steps differ from one since $\max_{k} \vert \lambda_k-1\vert> 0.4$. Nevertheless, the Newton algorithm still converge after $17$ iterates. More interestingly, the value $\beta=10^3$ illustrates the features and robustness of the algorithm: the convergence is achieved after 19 iterates. Far away from a zero of $E$, the variations of the error functional $E(y_k,f_k)$  are first quite slow, then increase to become very fast after 16 iterates, when $\lambda_k$ is close to one. In contrast, for $\beta=10^3.$, the Newton algorithm, still initialized with the linear solution diverges (see Table \ref{tab:LS_beta1000_N}). As discussed in \cite{lemoinemunch_time}, in that case, a continuation method with respect to the parameter $\beta$ may be combined with the Newton algorithm.  

On the contrary, we mention that with these data, the sequences obtained from the algorithm \eqref{NL_z_k} based on the linearization introduced in \cite{EFC-EZ}, remain bounded but do not converge, including for the value $\beta=10.$ The convergence is observed for instance with a larger size of the domain $\omega$, for instance $\omega=(0.2,0.8)$ (see \cite[section 4.2]{EFC-AM}).

\begin{table}[http]
	\centering
		\begin{tabular}{|c|c|c|c|c|c|c|}
			\hline
			$\sharp$iterate $k$  			   & $\frac{\Vert y_{k}-y_{k-1}\Vert_{L^2(Q_T)}}{\Vert y_{k-1}\Vert_{L^2(Q_T)}}$ & $\frac{\Vert f_{k}-f_{k-1}\Vert_{L^2(q_T)}}{\Vert f_{k-1}\Vert_{L^2(q_T)}}$ & $\Vert y_k\Vert_{2}$ & $\Vert f_k\Vert_{2,q_T}$ &  $\sqrt{2 E(y_k)}$   & $\lambda_k$ \tabularnewline
			\hline
			$0$ & $-$ & $-$ & $4.528$ & $4.391$ & $5.58\times 10^{-1}$ & $0.961$\tabularnewline
$1$ & $1.83\times 10^{-2}$ & $1.28\times 10^{-3}$ & $4.651$ & $4.402$ & $1.81\times 10^{-3}$ & $0.996$\tabularnewline
$2$ & $4.45\times 10^{-4}$ & $9.07\times 10^{-5}$ & $4.661$ & $4.403$ & $2.72\times 10^{-6}$ & $1.$\tabularnewline
$3$ & $1.12\times 10^{-6}$ & $3.74\times 10^{-7}$ & $4.662$ & $4.404$ & $4.88\times 10^{-8}$ & $1.$\tabularnewline

																				\hline
		\end{tabular}
	\caption{$\beta=10.$ ; Results for the algorithm (\ref{algo_LS_Y}).}
	\label{tab:LS_beta10}
\end{table}

\begin{table}[http]
	\centering
		\begin{tabular}{|c|c|c|c|c|c|c|}
			\hline
			$\sharp$iterate $k$  			   & $\frac{\Vert y_{k}-y_{k-1}\Vert_{L^2(Q_T)}}{\Vert y_{k-1}\Vert_{L^2(Q_T)}}$ & $\frac{\Vert f_{k}-f_{k-1}\Vert_{L^2(q_T)}}{\Vert f_{k-1}\Vert_{L^2(q_T)}}$ & $\Vert y_k\Vert_{2}$ & $\Vert f_k\Vert_{2,q_T}$ &  $\sqrt{2 E(y_k)}$   & $\lambda_k$ \tabularnewline
			\hline
			$0$ & $-$ & $-$ & $45.28$ & $43.91$ & $9.31\times 10^{-1}$ & $0.534$\tabularnewline
$1$ & $8.41\times 10^{-1}$ & $1.23\times 10^{-2}$ & $35.8908$ & $38.76$ & $1.12\times 10^{-1}$ & $0.591$\tabularnewline
$2$ & $1.93\times 10^{-1}$ & $2.91\times 10^{-3}$ & $36.7302$ & $38.92$ & $3.40\times 10^{-2}$ & $0.701$\tabularnewline
$3$ & $3.65\times 10^{-2}$ & $1.01\times 10^{-3}$ & $37.0919$ & $39.12$ & $6.12\times 10^{-3}$ & $0.812$\tabularnewline
$4$ & $1.12\times 10^{-2}$ & $2.69\times 10^{-4}$ & $37.2124$ & $40.01$ & $1.12\times 10^{-3}$ & $0.881$\tabularnewline
$5$ & $3.23\times 10^{-4}$ & $4.23\times 10^{-5}$ & $37.2426$ & $40.04$ & $2.13\times 10^{-4}$ & $0.912$\tabularnewline
$6$ & $1.27\times 10^{-5}$ & $6.23\times 10^{-6}$ & $37.2518$ & $40.05$ & $3.05\times 10^{-5}$ & $0.999$\tabularnewline
$7$ & $5.09\times 10^{-6}$ & $8.12\times 10^{-7}$ & $37.2520$ & $40.05$ & $2.10\times 10^{-6}$ & $0.999$\tabularnewline
$8$ & $7.40\times 10^{-8}$ & $8.21\times 10^{-9}$ & $37.2520$ & $40.05$ & $5.10\times 10^{-9}$ & $1.$\tabularnewline
																\hline
		\end{tabular}
	\caption{$\beta=10^2$ ; Results for the algorithm (\ref{algo_LS_Y}).}
	\label{tab:LS_beta100}
\end{table}

\begin{table}[http]
	\centering
		\begin{tabular}{|c|c|c|c|c|c|c|}
			\hline
				$\sharp$iterate $k$  			   & $\frac{\Vert y_{k}-y_{k-1}\Vert_{L^2(Q_T)}}{\Vert y_{k-1}\Vert_{L^2(Q_T)}}$ & $\frac{\Vert f_{k}-f_{k-1}\Vert_{L^2(q_T)}}{\Vert f_{k-1}\Vert_{L^2(q_T)}}$ & $\Vert y_k\Vert_{2}$ & $\Vert f_k\Vert_{2,q_T}$ &  $\sqrt{2 E(y_k)}$   & $\lambda_k$ \tabularnewline
			\hline
			$0$ 	    &   $-$	& $-$	                      &  $452.80$ & $439.18$ & $9.809\times10^{-1}$               & $0.4215$   \tabularnewline
			$1$       & $8.21\times 10^{-1}$ & $6.00\times 10^{-1}$ & $320.12$ & $330.15$ & $8.536\times10^{-1}$              & $0.3919$   \tabularnewline
			$2$ 	     &  $6.19\times 10^{-1}$ & $3.29\times 10^{-2}$ &324.02 & $ 334.12$ &  $8.012\times10^{-1}$ 	      &  $0.1566$  \tabularnewline
			$3$ 	     &  $4.18\times 10^{-1}$  & $1.37\times 10^{-2}$ & $325.65$ & $338.21$ & $7.953\times10^{-1}$ 		& $0.1767$   \tabularnewline
			$4$ 	     & $3.11\times 10^{-2}$  & $1.34\times 10^{-2}$ &  $326.11$ &  $340.12$ & $7.851\times10^{-1}$ 		 & $0.0937$   	\tabularnewline
			$5$ 	     &$2.98\times 10^{-2}$  & $5.85\times 10^{-3}$ & $326.35$ & $342.24$    & $7.688\times10^{-2}$ 		 & $0.0491$   	\tabularnewline
			$6$ 	    &  $3.37\times 10^{-2}$  & $7.00\times 10^{-3}$ & $326.91$ &  $344.65$  & $7.417\times10^{-2}$ 		 & $0.1296$   	\tabularnewline
			$7$ 	    &  $4.17\times 10^{-2}$ & $9.69\times 10^{-3}$ & $327.23$ & $346.12$  & $6.864\times10^{-2}$ 		 & $0.1077$   	\tabularnewline
			$8$ 	    & $2.89\times 10^{-2}$   & $8.09\times 10^{-3}$ & $327.42$ &  $347.19$ & $6.465\times10^{-2}$ 		 & $0.0859$   	\tabularnewline
			$9$ 	    & $1.09\times 10^{-2}$  & $6.40\times 10^{-3}$ &  $327.49$ & 347.29 & $6.182\times10^{-2}$ 		 & $0.0968$   \tabularnewline
			$10$ 	&    $1.02\times 10^{-2}$    & $6.72\times 10^{-3}$ & $327.92$ & 347.38  & $5.805\times10^{-2}$ 		 & $0.1184$   	\tabularnewline
			$11$ 	  &  $6.32\times 10^{-3}$    & $6.91\times 10^{-3}$ & $328.13$ & 347.41 & $5.371\times10^{-2}$ 		 & $0.1730$   	\tabularnewline
			$12$ 	   & $5.53\times 10^{-3}$    & $7.41\times 10^{-3}$ & $328.16$ & 347.43 & $4.825\times10^{-2}$ 		 & $0.2579$   	\tabularnewline
			$13$ 	    & $4.32\times 10^{-3}$   & $8.22\times 10^{-3}$ & $328.19$  & 347.45& $4.083\times10^{-2}$ 		 & $0.3817$   	\tabularnewline
			$14$ 	  &  $2.13\times 10^{-3}$    & $8.14\times 10^{-3}$ & $328.21$ & 347.48 & $3.164\times10^{-2}$ 		 & $0.4946$   	\tabularnewline
			$15$ 	  &  $3.57\times 10^{-3}$    & $7.34\times 10^{-3}$ &$328.22 $ & 347.50  & $2.207\times10^{-2}$ 		 & $0.8294$   	\tabularnewline
			$16$ 	  &  $1.01\times 10^{-3}$    & $6.68\times 10^{-3}$ &$328.25$ & 347.51 & $1.174\times10^{-2}$ 		 & $0.9845$   	\tabularnewline
			$17$ 	  &  $5.68\times 10^{-4}$    & $3.84\times 10^{-4}$ & $328.26$ & 347.51 & $2.191\times10^{-3}$ 		 & $0.9999$  	\tabularnewline
			$18$ 	  &   $2.14\times 10^{-4}$   & $5.85\times 10^{-5}$ &$328.26$& 347.52 & $4.674\times10^{-5}$ 		 & $1.$   	\tabularnewline
			$19$ 	  &   $3.21\times 10^{-6}$   & $1.57\times 10^{-7}$ &$328.27$ & 347.52 & $5.843\times10^{-7}$ 		 & $-$  	\tabularnewline
									\hline
		\end{tabular}
	\caption{$\beta=10^3$ ; Results for the algorithm (\ref{algo_LS_Y}).}
	\label{tab:LS_beta1000}

\end{table}

\begin{table}[http]
	\centering
		\begin{tabular}{|c|c|c|c|c|c|}
			\hline
				$\sharp$iterate $k$  			   & $\frac{\Vert y_{k}-y_{k-1}\Vert_{L^2(Q_T)}}{\Vert y_{k-1}\Vert_{L^2(Q_T)}}$ & $\frac{\Vert f_{k}-f_{k-1}\Vert_{L^2(q_T)}}{\Vert f_{k-1}\Vert_{L^2(q_T)}}$ & $\Vert y_k\Vert_{2}$ & $\Vert f_k\Vert_{2,q_T}$ &  $\sqrt{2 E(y_k)}$   \tabularnewline
			\hline
			$0$ 	    &   $-$	& $-$	                      &  $452.80$ & $439.18$ & $9.809\times10^{-1}$                  \tabularnewline
			$1$       & $9.76\times 10^{-1}$ & $1.05$ & $330.21$ & $334.15$ & $9.812\times10^{-1}$                \tabularnewline
			$2$ 	     &  $1.02$ & $1.11$ &344.37 & $ 336.12$ &  $1.356$ 	       \tabularnewline
			$3$ 	     &  $1.27$  & $1.13$ & $366.92$ & $338.23$ & $4.319$ 	  \tabularnewline
			$4$ 	     & $1.18$  & $1.25$ &  $406.06$ &  $343.12$ & $4.799$ 		  	\tabularnewline
			$5$ 	     &$1.01$  & $1.14$ & $481.53$ & $405.03$    & $13.131$ 		   	\tabularnewline
		
												\hline
		\end{tabular}
	\caption{$\beta=10^3$ ; Results for the algorithm (\ref{algo_LS_Y}) with $\lambda_k=1$ for all $k$.}
	\label{tab:LS_beta1000_N}

\end{table}

\section{Conclusions and perspectives}\label{conclusion}

We have constructed an explicit sequence of functions $(f_k)_{k}$ converging strongly in the $L^2(q_T)$ norm toward a null control for the semilinear heat equation $y_t-\Delta y+ g(y)=f \, 1_\omega$. The construction of the sequence is based on the minimization of a $L^2(0,T;H^{-1}(\Omega))$ least-squares functional. The use of a specific descent direction allows to achieve a global convergence (uniform with respect to the data and to the initial guess) with a super-linear rate related to the regularity of the nonlinear function $g$. Experiment confirms the robustness of the approach. In this analysis, we have assumed in particular that the derivative $g^{\prime}$ of $g$ is uniformly bounded in $\mathbb{R}$. This allows to get a uniform bound of the constant of the form $C(\Omega,\omega,T,\Vert g^{\prime}(y)\Vert_{\infty})$ appearing  from the Carleman estimate \eqref{Carleman-ine}. 
In order to remove this assumption and be able to consider super-linear function $g$ (as in the seminal work \cite{EFC-EZ} by Fern\'andez-Cara and Zuazua, assuming that $g$ is locally Lipschitz-continuous and the asymptotic behavior \eqref{asymptotic_g}), we need to refine the analysis and exploit the structure of the constant $C(\Omega,\omega,T,\Vert g^{\prime}(y)\Vert_{L^\infty})$ (as done in \cite{Duyckaerts} for the observability constant).  This may allow, assuming the above hypotheses of \cite{EFC-EZ}, not only to recover the null controllability of \eqref{heat-NL} but also, to construct, within the algorithm (\ref{algo_LS_Y}), approximations of null controls. 

We also emphasize that this least-squares approach is very general and may be used to address other PDEs. Following \cite{lemoinemunch_time} devoted the direct problem, one may notably study the applicability of the method to approximate control for the Navier-Stokes system. We also mentioned the case of nonlinear wave equation studied in \cite{Zuazua_WaveNL} making use of a fixed point strategy.

 \bibliographystyle{plain}

\end{document}